\documentclass[11pt,a4paper,draft]{article}
\usepackage{fullpage}
\usepackage{amssymb}
\usepackage{amsmath}
\usepackage{amsthm}
\usepackage{amsrefs}
\usepackage[arrow,matrix,curve,tips]{xy}
\SelectTips{cm}{}
\usepackage{tikz}
\let\mod=\undefined
\DeclareMathOperator{\AR}{AR}
\DeclareMathOperator{\dist}{dist}
\DeclareMathOperator{\dimv}{\mathbf{dim}}
\DeclareMathOperator{\End}{End}
\DeclareMathOperator{\Ext}{Ext}
\DeclareMathOperator{\GL}{GL}
\DeclareMathOperator{\Hom}{Hom}

\DeclareMathOperator{\mod}{mod}
\DeclareMathOperator{\mult}{mult}
\DeclareMathOperator{\op}{op}
\DeclareMathOperator{\rk}{rk}
\DeclareMathOperator{\rep}{rep}
\DeclareMathOperator{\Spec}{Spec}
\DeclareMathOperator{\supp}{supp}
\newcommand{\BA}{{\mathbb A}}
\newcommand{\BB}{{\mathbb B}}

\newcommand{\BD}{{\mathbb D}}
\newcommand{\BE}{{\mathbb E}}
\newcommand{\BM}{{\mathbb M}}
\newcommand{\BN}{{\mathbb N}}

\newcommand{\BZ}{{\mathbb Z}}
\newcommand{\CA}{{\mathcal A}}

\newcommand{\CC}{{\mathcal C}}
\newcommand{\CD}{{\mathcal D}}

\newcommand{\CI}{{\mathcal I}}
\newcommand{\CL}{{\mathcal L}}
\newcommand{\CM}{{\mathcal M}}

\newcommand{\CO}{{\mathcal O}}
\newcommand{\CP}{{\mathcal P}}

\newcommand{\CR}{{\mathcal R}}
\newcommand{\CS}{{\mathcal S}}

\newcommand{\CV}{{\mathcal V}}

\newcommand{\CX}{{\mathcal X}}
\newcommand{\CY}{{\mathcal Y}}

\newcommand{\dd}{{\mathbf d}}
\newcommand{\ee}{{\mathbf e}}
\newcommand{\hh}{{\mathbf h}}
\newcommand{\mm}{{\mathbf m}}
\newcommand{\vv}{{\mathbf v}}
\newcommand{\frakb}{{\mathfrak b}}
\newcommand{\frakq}{{\mathfrak q}}
\newcommand{\ov}{\overline}
\newcommand{\un}{\underline}
\newcommand{\bsmatrix}[1]{\left[\begin{smallmatrix} #1%
 \end{smallmatrix}\right]}
\newtheorem{thm}{Theorem}[section]
\newtheorem{cor}[thm]{Corollary}
\newtheorem{lem}[thm]{Lemma}
\newtheorem{prop}[thm]{Proposition}
\theoremstyle{definition}
\newtheorem{ex}[thm]{Example}
\newtheorem{df}[thm]{Definition}
\numberwithin{equation}{section}
\begin{document}

\title{Tangent spaces of orbit closures for representations of Dynkin quivers of type $\BD$
}
\author{Grzegorz Bobi\'nski and Grzegorz Zwara}
\date{}

\maketitle

\begin{abstract}
Let $\Bbbk$ be an algebraically closed field, $Q$ a finite quiver,
and denote by $\rep_Q^\dd$ the affine $\Bbbk$-scheme of representations of $Q$
with a fixed dimension vector $\dd$.
Given a representation $M$ of $Q$ with dimension vector $\dd$, the set
$\CO_M$ of points in $\rep_Q^\dd(\Bbbk)$ isomorphic as representations to $M$
is an orbit under an action on $\rep_Q^\dd(\Bbbk)$ of a product of general linear groups.
The orbit $\CO_M$ and its Zariski closure $\ov{\CO}_M$, considered as reduced subschemes of $\rep_Q^\dd$,
are contained in an affine scheme $\CC_M$ defined
by rank conditions on suitable matrices associated to $\rep_Q^\dd$.
For all Dynkin and extended Dynkin quivers, the sets of points of $\ov{\CO}_M$ and $\CC_M$ coincide,
or equivalently, $\ov{\CO}_M$ is the reduced scheme associated to $\CC_M$.
Moreover, $\ov{\CO}_M=\CC_M$ provided $Q$ is a Dynkin quiver of type $\BA$,
and this equality is a conjecture for the remaining Dynkin quivers (of type $\BD$ and $\BE$).

Let $Q$ be a Dynkin quiver of type $\BD$ and $M$ a finite dimensional representation of $Q$.
We show that the equality $T_N\ov{\CO}_M=T_N\CC_M$ of Zariski tangent spaces holds
for any closed point $N$ of $\ov{\CO}_M$.
As a consequence, we describe the tangent spaces to $\ov{\CO}_M$ in representation theoretic terms.
\end{abstract}

\section{Introduction and the main results}
\label{intro}

Throughout the paper $\Bbbk$ denotes an algebraically closed field
of arbitrary characteristic.
We will identify a $\Bbbk$-scheme $\CX$ with its functor of points,
i.e.\ the functor from the category of commutative $\Bbbk$-algebras
to the category of sets sending $R$ to the set of morphisms
$\Spec(R)\to\CX$.
Let $\Bbbk[\varepsilon]$ denote the $\Bbbk$-algebra of dual numbers and consider the map
\[
\CX(\pi)\colon\CX(\Bbbk[\varepsilon])\to\CX(\Bbbk),
\]
where $\pi\colon\Bbbk[\varepsilon]\to\Bbbk$ is the canonical surjective homomorphism.
Given a $\Bbbk$-rational point $x$ of $\CX$, i.e.\ $x\in\CX(\Bbbk)$,
the fiber $\CX(\pi)^{-1}(x)$ is the Zariski tangent space $T_x\CX$ to $\CX$ at $x$.
We are mostly interested in affine schemes $\CX$ of finite type over $\Bbbk$, i.e.\
the schemes of the form $\Spec(R)$, where $R$ is a finitely generated commutative $\Bbbk$-algebra.
For such schemes $\CX(\Bbbk)$ is the set of closed points of $\CX$.

Let $Q=(Q_0,Q_1)$ be a finite quiver, i.e.\ a finite set $Q_0$
of vertices and a finite set $Q_1$ of arrows $\alpha\colon s\alpha\to t\alpha$,
where $s\alpha$ and $t\alpha$ denote the starting
and the terminating vertex of $\alpha$, respectively.
A representation of $Q$ over $\Bbbk$ is a collection $M=(M_a,M_\alpha;\,a\in Q_0,\alpha\in Q_1)$
of $\Bbbk$-vector spaces $M_a$ and $\Bbbk$-linear maps $M_\alpha\colon M_{s\alpha}\to M_{t\alpha}$.
A morphism $f\colon M\to N$ between two representations is
a collection $f=(f_a\colon M_a\to N_a;\,a\in Q_0)$ of $\Bbbk$-linear maps such that
\begin{equation*}
f_{t\alpha}\circ M_\alpha=N_\alpha\circ f_{s\alpha},\qquad
\text{for all $\alpha\in Q_1$}.
\end{equation*}
We denote by $\rep(Q)$ the category of finite dimensional representations of $Q$,
i.e.\ the representations $M$ such that all vector spaces $M_a$ are finite dimensional.
For a representation $M$ in $\rep(Q)$ we define its dimension vector $\dimv M=(\dim_\Bbbk M_a)\in\BN^{Q_0}$.

We denote by $\BM_{p,q}$ the $\Bbbk$-scheme of $p\times q$-matrices and by $\GL_d$ the group $\Bbbk$-scheme
of invertible $d\times d$-matrices, for any positive integers $p$, $q$ and $d$. Given a dimension vector $\dd=(\dd_a)\in\BN^{Q_0}$ we have the affine scheme
\[
\rep_Q^\dd=\prod_{\alpha\in Q_1}\BM_{\dd_{t\alpha}\times \dd_{s\alpha}}.
\]
Thus the points of $\rep_Q^\dd(\Bbbk)$ can be identified with the representations $M$ of $Q$ such that $M_a=\Bbbk^{\dd_a}$
for any $a\in Q_0$.
The group scheme
\[
\GL_\dd=\prod_{a\in Q_0}\GL_{\dd_a}
\]
acts on $\rep_Q^\dd$ via
\[
g\star M=(g_{t\alpha}\cdot M_\alpha\cdot g_{s\alpha}^{-1}),
\]
for any $g=(g_a)\in\GL_\dd(R)$, $M=(M_\alpha)\in\rep_Q(R)$, and a commutative $\Bbbk$-algebra $R$.
Given a representation $M$ in $\rep(Q)$, we denote by $\CO_M$ the $\GL_\dd(\Bbbk)$-orbit in $\rep_Q^\dd(\Bbbk)$
which consists of the representations in $\rep_Q^\dd(\Bbbk)$ isomorphic to $M$, where $\dd=\dimv M$.
By abuse of notation, we treat $\CO_M$ and its closure $\ov{\CO}_M$
as reduced subschemes of $\rep_Q^\dd$.
It is an open and important problem to describe the defining ideal of $\ov{\CO}_M$ or even to exhibit polynomials
having $\ov{\CO}_M$ as their zero set.
If $M$ and $N$ are representations satisfying $\CO_N\subseteq\ov{\CO}_M$, then
we say that $M$ degenerates to $N$.
Note that $\ov{\CO}_M (\Bbbk)$ is the union of $\CO_N (\Bbbk)$,
where $N$ runs through the representations to which $M$ degenerates.

Let $[X,Y]=\dim_\Bbbk\Hom_Q(X,Y)$, for $X,Y\in\rep(Q)$.
It is well-known that if $M$ degenerates to $N$ then
\begin{equation} \label{homorder}
[X,N]\geq[X,M]\quad\text{and}\quad[N,X]\geq[M,X],\qquad\text{for any $X\in\rep(Q)$}
\end{equation}
(see for instance~\cite{Rie}).
Moreover, due to Bongartz~\cites{Bdeg, Bext}, the reverse implication holds under an additional assumption
on $Q$.

\begin{thm} \label{degeqhom}
Let $Q$ be a Dynkin or an extended Dynkin quiver.
Assume $M$ and $N$ belong to $\rep(Q)$ and $\dimv M=\dimv N$.
Then $M$ degenerates to $N$ if and only if the condition \eqref{homorder} is satisfied.
\qed
\end{thm}

Inspired by the above inequalities (see also~\cite{Bdeg}*{Proposition 1}),
a closed $\GL_\dd$-subscheme $\CC_M$ of $\rep_Q^\dd$ containing $\ov{\CO}_M$
was defined in~\cite{RZrank}.
Let $\Bbbk Q=\bigoplus_{a,b\in Q_0}\Bbbk Q(a,b)$ denote the path algebra of $Q$,
where $\Bbbk Q(a,b)$ is the vector space with a $\Bbbk$-basis formed by the paths in $Q$
starting at $b$ and terminating at $a$.
For any commutative $\Bbbk$-algebra $R$, $X\in\rep_Q^\dd(R)$ and $\omega\in\Bbbk Q(a,b)$,
the matrix $X_\omega\in\BM_{\dd_a,\dd_b}(R)$ is defined in the obvious way.

Let $p,q\in\BN$, and consider two sequences $(a_1,\ldots,a_p)$ and
$(b_1,\ldots,b_q)$ of vertices in $Q_0$ and a $p\times q$-matrix
$\un{\omega}=(\omega_{i,j})$ such that each $\omega_{i,j}$ belongs to $\Bbbk Q(a_i,b_j)$.
We define a regular morphism
\[
\Theta_{\un{\omega}}\colon\rep_Q^\dd \to \BM_{p' \times q'},\qquad
\Theta_{\un{\omega}}(N)=\begin{bmatrix}
N_{\omega_{1,1}}&\cdots&N_{\omega_{1,q}}\\
\hdotsfor{3}\\
N_{\omega_{p,1}}&\cdots&N_{\omega_{p,q}}
\end{bmatrix},
\]
where $p'=\sum\dd_{a_i}$ and $q'=\sum\dd_{b_j}$.
For $M$ in $\rep(Q)$, let $\dd=\dimv M$ and $\CI_{M,\un{\omega}}$ be the ideal in the coordinate algebra $\Bbbk[\rep_Q^\dd]$
of $\rep_Q^\dd$ generated by the images via $(\Theta_{\un{\omega}})^\ast$ of the minors of size $1+\rk\Theta_{\un{\omega}}(M')$
in $\Bbbk[\BM_{p' \times q'}]$,
where $M'$ is a point of $\CO_M$.
We set $\CI_M=\sum\CI_{M,\un{\omega}}$, where $\un{\omega}$ runs through
all possible matrices
of linear combinations of paths with all possible sequences of starting and terminating vertices.
Then $\CC_M=\Spec(\Bbbk[\rep_Q^\dd]/\CI_M)$ is a closed $\GL_\dd$-subscheme of $\rep_Q^\dd$ containing $\ov{\CO}_M$.
Moreover, $\CC_M(\Bbbk)$ consists of the representations $N\in\rep_Q^\dd(\Bbbk)$ satisfying \eqref{homorder}.
Hence we can reformulate Theorem~\ref{degeqhom} by saying that $\CC_M(\Bbbk)=\ov{\CO}_M(\Bbbk)$
provided $Q$ is a Dynkin or an extended Dynkin quivers.
Note that this equality does not hold for the representation
\[
M=\quad
\xymatrix{\Bbbk\ar@<.5ex>[r]^{\bsmatrix{1\\0}}\ar@<-.5ex>[r]_{\bsmatrix{0\\1}}
 &\Bbbk^2\ar@<.5ex>[r]^{\bsmatrix{1&0}} \ar@<-.5ex>[r]_{\bsmatrix{0&1}}&\Bbbk.}
\]
Here $\CC_M(\Bbbk)$ has two irreducible components of dimension $5$, with one of them being $\ov{\CO}_M(\Bbbk)$
(see~\cite{RZrank}*{Example 8.9}).

Lakshmibai and Magyar described in~\cite{LM} generators of the defining ideal of $\ov{\CO}_M$ in case $Q$ is an equioriented
Dynkin quiver of type $\BA$.
They used the Zelevinsky immersion of $\rep_Q^\dd$ in a Schubert variety of a flag variety,
and applied a description of generators of the defining ideal of this Schubert variety.
It turned out that they showed that the defining ideal of $\ov{\CO}_M$ equals $\CI_M$.
The result was generalized in~\cite{RZrank} to the Dynkin quivers of type $\BA$ with an arbitrary orientation:

\begin{thm} \label{mainforA}
Let $Q$ be a Dynkin quiver of type $\BA$ and $M\in\rep(Q)$.
Then $\CC_M=\ov{\CO}_M$.
\qed
\end{thm}

We do not know if this result remains true for the Dynkin quivers
of type $\BD$ ($\BD_n$, $n\geq 4$) and $\BE$ ($\BE_6$, $\BE_7$, $\BE_8$).
It would be interesting to know if $\CC_M$ and $\ov{\CO}_M$ have at least identical
tangent spaces in these cases.
Note that this is not true even for the simplest extended Dynkin quiver.
Namely, consider representations
\[
M=\quad\xymatrix@1{\Bbbk^2\ar@(ur,dr)[]^{\bsmatrix{0&1\\ 0&0}}}
\qquad \text{and} \qquad
N=\quad\xymatrix@1{\Bbbk^2\ar@(ur,dr)[]^{\bsmatrix{0&0\\ 0&0}}}.
\]
Then $N\in\ov{\CO}_M$, $T_N\ov{\CO}_M$ has dimension $3$, while $T_N\CC_M$ has dimension $4$
(the trace function does not belong to the ideal $\CI_M$, see~\cite{RZrank}*{Example 3.7}).

Our main result shows that $\CC_M$ and $\ov{\CO}_M$ have identical tangent spaces
in type $\BD$:

\begin{thm} \label{main}
Let $Q$ be a Dynkin quiver of type $\BD$ and $M\in\rep(Q)$.
Then $\ov{\CO}_M(\Bbbk[\varepsilon])=\CC_M(\Bbbk[\varepsilon])$.
In other words, $T_N\ov{\CO}_M=T_N\CC_M$ for any $N$ in $\ov{\CO}_M(\Bbbk)=\CC_M(\Bbbk)$.
\end{thm}

We give now a representation theoretic interpretation of $T_N\CC_M$.
The sets $\rep_Q^\dd(\Bbbk)$ and $\rep_Q^\dd(\Bbbk[\varepsilon])$ have natural structures
of vector spaces over $\Bbbk$.
Using the decomposition $\Bbbk[\varepsilon]=\Bbbk\cdot 1\oplus\Bbbk\cdot\varepsilon$
one can write each element of $\rep_Q^\dd(\Bbbk[\varepsilon])$ uniquely in the form
$N+\varepsilon\cdot Z$, where $N$ and $Z$ belong to $\rep_Q^\dd(\Bbbk)$.
Here $N$ is a closed point of $\rep_Q^\dd$ and $Z$ is a tangent vector in $T_N\rep_Q^\dd$.
It is a simple but ingenious idea to associate to such element $N+\varepsilon\cdot Z$ a representation
in $\rep_Q^{2\cdot\dd}(\Bbbk)$ of the following block form
\[
\bsmatrix{N&Z\\ 0&N}=\left(\bsmatrix{N_\alpha&Z_\alpha\\ 0&N_\alpha};\,\alpha\in Q_1\right).
\]
By~\cite{RZrank}, if $N\in\CC_M(\Bbbk)$, then $T_N\CC_M$ consists of the $Z$ satisfying
the following two equivalent conditions:
\begin{itemize}
\item $\bigl[X,\bsmatrix{N&Z\\ 0&N}\bigr]=2\cdot[X,N]$ for all $X\in\rep(Q)$ with $[X,N]=[X,M]$,
\item $\bigl[\bsmatrix{N&Z\\ 0&N},X\bigr]=2\cdot[N,X]$ for all $X\in\rep(Q)$ with $[N,X]=[M,X]$.
\end{itemize}

We remark that Theorem~\ref{main} should have applications in the problem of describing
the singular locus of $\ov{\CO}_M$ in representation theoretic terms,
for the representations $M$ of the Dynkin quivers of type $\BD$ (see~\cite{RZrank}*{Section 8}).

The paper is organized as follows:
in Section~\ref{sect two} we prove necessary facts about exact sequences in $\rep (Q)$
(in fact we formulate them in a more general setup of triangles in the derived category of $\rep (Q)$),
while in Section~\ref{sectionproof} we apply results of Section~\ref{sect two}
in geometric context and prove the main result.
For basic background on representation theory of quivers we refer to~\cites{ARS, ASS, Rin}.

The both authors gratefully acknowledge the support of the National Science Centre grant no.~2020/37/B/ST1/00127.

\section{Derived categories for representations of Dynkin quivers}
\label{sect two}

In order to prove Theorem~\ref{main}, we need a result about existence of short exact sequences
in $\rep(Q)$ with some special properties, where $Q$ is a Dynkin quiver of type $\BD$ (Corollary~\ref{maincortriang}).
Our idea is to use the embedding of $\rep(Q)$ in its derived category $\CD^b(Q)=\CD^b(\rep(Q))$,
and to prove existence of triangles in $\CD^b(Q)$ satisfying similar properties (Proposition~\ref{mainproptriang}).
An advantage of working with the derived category is that its structure, including Auslander-Reiten theory,
is more ``regular'' than the structure of $\rep(Q)$.
We refer to~\cite{Hap} as a general reference for this section.

\subsection{Dynkin graphs}

Throughout this section $\Delta=(\Delta_0,\Delta_1)$ is a Dynkin graph
of one of the types $\BA$, $\BD$ or $\BE$:
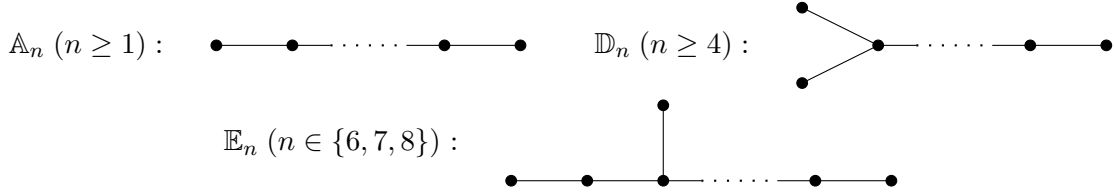
\begin{gather*}
\BA_n\;(n\geq 1):
\quad
\vcenter{\hbox{
\begin{tikzpicture}
\draw (0,0) -- (1.5,0);
\draw [loosely dotted,thick] (1.5,0) -- (2.5,0);
\draw (2.5,0) -- (4,0);
\filldraw
 (0,0) circle (2pt)
 (1,0) circle (2pt)
 (3,0) circle (2pt)
 (4,0) circle (2pt);
\end{tikzpicture}
}}
\qquad
\BD_n\;(n\geq 4):
\quad
\vcenter{\hbox{
\begin{tikzpicture}
\draw (0,.5) -- (1,0);
\draw (0,-.5) -- (1,0);
\draw (1,0) -- (1.5,0);
\draw [loosely dotted,thick] (1.5,0) -- (2.5,0);
\draw (2.5,0) -- (4,0);
\filldraw
 (0,.5) circle (2pt)
 (0,-.5) circle (2pt)
 (1,0) circle (2pt)
 (3,0) circle (2pt)
 (4,0) circle (2pt);
\end{tikzpicture}
}} \\
\BE_n\;(n\in\{6,7,8\}):
\quad
\vcenter{\hbox{
\begin{tikzpicture}
\draw (0,0) -- (2.5,0);
\draw (2,0) -- (2,1);
\draw [loosely dotted,thick] (2.5,0) -- (3.5,0);
\draw (3.5,0) -- (5,0);
\filldraw
 (0,0) circle (2pt)
 (1,0) circle (2pt)
 (2,0) circle (2pt)
 (2,1) circle (2pt)
 (4,0) circle (2pt)
 (5,0) circle (2pt);
\end{tikzpicture}
}}
\end{gather*}
where $\Delta_0$ is the set of $n$ vertices of $\Delta$, and $\Delta_1$ is its set of edges,
i.e.\ two element subsets of $\Delta_0$.
If $\{a,b\}$ is an edge we say that $a$ and $b$ are adjacent.
We denote by $a^-$ the (open) neighbourhood of a vertex $a$,
i.e.\ the set of vertices adjacent to $a$.
The degree of $a$ equals, by definition, the cardinality of $a^-$.
Let $\dist(a,b)$ be the length of the shortest walk in $\Delta$ between $a$ and $b$.

We define an integer $n_\Delta$ as follows:
\begin{equation} \label{nDelta}
n_{\BA_n}=n+1,\quad
n_{\BD_n}=2n-2,\quad
n_{\BE_6}=12,\quad
n_{\BE_7}=18,\quad
n_{\BE_8}=30.
\end{equation}
With $\Delta$ we associate (Tits) quadratic form
\begin{equation} \label{qDelta}
\frakq_\Delta\colon\BZ^{\Delta_0}\to\BZ,\quad
\frakq_\Delta(\dd)=\sum_{a\in\Delta_0}\dd_a^2-\sum_{\{a,b\}\in\Delta_1}\dd_a\cdot\dd_b,
\end{equation}
which is positive definite, i.e.\ $\frakq_\Delta(\dd)>0$ for any non-zero $\dd$.
If $\frakq_\Delta(\dd)=1$ we say $\dd$ is a root.
Obviously, if $\dd$ is a root, then $-\dd$ is also a root.
There are $n\cdot n_\Delta$ roots, half of them are positive,
where a vector $\dd$ is called positive provided $\dd\neq 0$ and $\dd_a\geq 0$, for each $a$.
There is a unique maximal root $\hh^\Delta$ (i.e.\ $\hh^\Delta-\dd$ is positive for any root $\dd\neq\hh^\Delta$) which equals
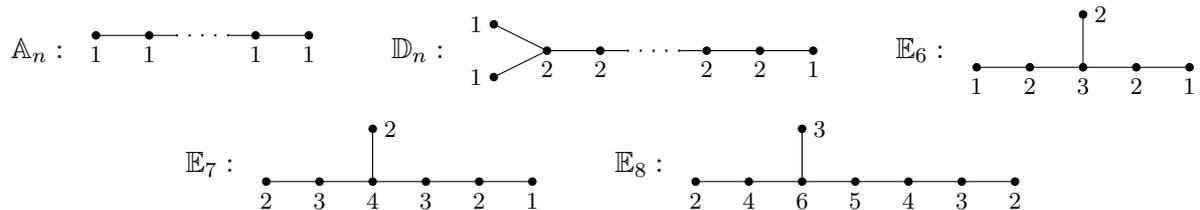
\begin{gather*}
\BA_n:
\vcenter{\hbox{
\begin{tikzpicture}[scale=.7] \footnotesize
\draw (0,0) -- (1.5,0);
\draw [loosely dotted,thick] (1.5,0) -- (2.5,0);
\draw (2.5,0) -- (4,0);
\filldraw
 (0,0) circle (2pt) node [below=1pt] {$1$}
 (1,0) circle (2pt) node [below=1pt] {$1$}
 (3,0) circle (2pt) node [below=1pt] {$1$}
 (4,0) circle (2pt) node [below=1pt] {$1$};
\end{tikzpicture}
}}
\qquad
\BD_n:
\vcenter{\hbox{
\begin{tikzpicture}[scale=.7] \footnotesize
\draw (0,.5) -- (1,0);
\draw (0,-.5) -- (1,0);
\draw (1,0) -- (2.5,0);
\draw [loosely dotted,thick] (2.5,0) -- (3.5,0);
\draw (3.5,0) -- (6,0);
\filldraw
 (0,.5) circle (2pt) node [left=1pt] {$1$}
 (0,-.5) circle (2pt) node [left=1pt] {$1$}
 (1,0) circle (2pt) node [below=1pt] {$2$}
 (2,0) circle (2pt) node [below=1pt] {$2$}
 (4,0) circle (2pt) node [below=1pt] {$2$}
 (5,0) circle (2pt) node [below=1pt] {$2$}
 (6,0) circle (2pt) node [below=1pt] {$1$};
\end{tikzpicture}
}}
\qquad
\BE_6:
\vcenter{\hbox{
\begin{tikzpicture}[scale=.7] \footnotesize
\draw (0,0) -- (4,0);
\draw (2,0) -- (2,1);
\filldraw
 (0,0) circle (2pt) node [below=1pt] {$1$}
 (1,0) circle (2pt) node [below=1pt] {$2$}
 (2,0) circle (2pt) node [below=1pt] {$3$}
 (2,1) circle (2pt) node [right=1pt] {$2$}
 (3,0) circle (2pt) node [below=1pt] {$2$}
 (4,0) circle (2pt) node [below=1pt] {$1$};
\end{tikzpicture}
}} \\
\BE_7:
\vcenter{\hbox{
\begin{tikzpicture}[scale=.7] \footnotesize
\draw (0,0) -- (5,0);
\draw (2,0) -- (2,1);
\filldraw
 (0,0) circle (2pt) node [below=1pt] {$2$}
 (1,0) circle (2pt) node [below=1pt] {$3$}
 (2,0) circle (2pt) node [below=1pt] {$4$}
 (2,1) circle (2pt) node [right=1pt] {$2$}
 (3,0) circle (2pt) node [below=1pt] {$3$}
 (4,0) circle (2pt) node [below=1pt] {$2$}
 (5,0) circle (2pt) node [below=1pt] {$1$};
\end{tikzpicture}
}}
\qquad
\BE_8:
\vcenter{\hbox{
\begin{tikzpicture}[scale=.7] \footnotesize
\draw (0,0) -- (6,0);
\draw (2,0) -- (2,1);
\filldraw
 (0,0) circle (2pt) node [below=1pt] {$2$}
 (1,0) circle (2pt) node [below=1pt] {$4$}
 (2,0) circle (2pt) node [below=1pt] {$6$}
 (2,1) circle (2pt) node [right=1pt] {$3$}
 (3,0) circle (2pt) node [below=1pt] {$5$}
 (4,0) circle (2pt) node [below=1pt] {$4$}
 (5,0) circle (2pt) node [below=1pt] {$3$}
 (6,0) circle (2pt) node [below=1pt] {$2$};
\end{tikzpicture}
}}
\end{gather*}

\subsection{Derived category for acyclic quivers}

Throughout this subsection $Q$ is a finite quiver without oriented cycles.
We denote by $\CD^b(Q)=\CD^b(\rep(Q))$ the derived category of the abelian category $\rep(Q)$.
The category $\CD^b(Q)$ is triangulated, hence there is an auto-equivalence $[1]$ of $\CD^b(Q)$ called
the shift functor (``the suspension functor'' and ``the translation functor'' are alternative names used by other authors)
and a class of triangles (the name ``distinguished triangles'' is commonly used),
written in the form $A\xrightarrow{\alpha}B\xrightarrow{\beta}C\xrightarrow{\gamma}A[1]$.
There is a canonical full embedding of $\rep(Q)$ in $\CD^b(Q)$,
and we shall identify $\rep(Q)$ with its image in $\CD^b(Q)$.
In particular,
\[
\Hom_Q(X,Y)=\Hom_{\CD^b(Q)}(X,Y)
\quad\text{and}\quad
\Ext^1_Q(X,Y)=\Hom_{\CD^b(Q)}(X,Y[1]),
\]
for all $X,Y\in\rep(Q)$.
Based on the latter equality, there is a strong relationship between the short exact sequences
in $\rep(Q)$ and triangles in $\CD^b(Q)$.
Namely, for each short exact sequence $\sigma\colon0\to A\xrightarrow{\alpha}B\xrightarrow{\beta}C\to 0$ in $\rep(Q)$
there is a unique morphism $\gamma\in\Hom_{\CD^b(Q)}(C,A[1])$ such that
$\hat{\sigma}\colon A\xrightarrow{\alpha}B\xrightarrow{\beta}C\xrightarrow{\gamma}A[1]$ is a triangle in $\CD^b(Q)$.
Conversely, if $A\xrightarrow{\alpha}B\xrightarrow{\beta}C\xrightarrow{\gamma}A[1]$ is a triangle in $\CD^b(Q)$
with $A$ and $C$ in $\rep(Q)$, then there is an isomorphism $g\colon B\to B'$ in $\CD^b(Q)$ such that
$\sigma\colon0\to A\xrightarrow{g\circ\alpha}B'\xrightarrow{\beta\circ g^{-1}}C\to 0$  is a short exact sequence in $\rep(Q)$.
Moreover, in the above situation $\hat{\sigma}$ has the form
$A\xrightarrow{g\circ\alpha}B'\xrightarrow{\beta\circ g^{-1}}C\xrightarrow{\gamma}A[1]$.

We generalize now to triangles notion of a split exact sequence and a pullback.
Let $\sigma$ be a triangle $A\xrightarrow{\alpha}B\xrightarrow{\beta}C\xrightarrow{\gamma}A[1]$ in $\CD^b(Q)$.
We say that $\sigma$ splits if one of the following equivalent conditions is satisfied: (i)~$\alpha$ is a section, (ii)~$\beta$ is a retraction,
(iii)~$\gamma=0$, (iv)~$B$ is isomorphic to $A\oplus C$.
Observe that an exact sequence $\sigma$ splits if and only if the triangle $\hat{\sigma}$ splits.

Given a triangle $\sigma\colon A\xrightarrow{\alpha}B\xrightarrow{\beta}C\xrightarrow{\gamma}A[1]$
and a morphism $h\colon C'\to C$ in $\CD^b(Q)$ we get the following commutative diagram
\[
\xymatrixcolsep{1pc}
\xymatrix{
A\ar[rr]\ar@{=}[d]&&B'\ar[rr]\ar[d]&&C'\ar[rr]^-{\gamma\circ h}\ar[d]^-h&&A[1]\ar@{=}[d]\\
A\ar[rr]^-{\alpha}&&B\ar[rr]^-{\beta}&&C\ar[rr]^-{\gamma}&&A[1],
}
\]
where the upper row, called the pullback of $\sigma$ along $h$, is a triangle
(note that the pullback is unique up to isomorphism of triangles).
One defines pushouts dually.
Observe that if $\sigma'$ is the pullback of a short exact sequence $\sigma$ along a homomorphism $h\colon C'\to C$ in $\rep(Q)$,
then $\hat{\sigma}'$ is the pullback of $\hat{\sigma}$ along $h$ (viewed as a morphism in $\CD^b(Q)$).

We say that a triangle $\sigma$ is almost split (or an Auslander-Reiten triangle)
if $A$ and $C$ are indecomposable, $\sigma$ does not split, but its pullbacks split for all
morphisms to $C$ which are not retractions.
The last condition can be replaced by the requirement that the pushouts of $\sigma$ split for all
morphisms from $A$ which are not sections.

An important fact about the category $\CD^b(Q)$ is that it has a Serre duality, i.e.\
there is an auto-equivalence $\nu\colon\CD^b(Q)\to\CD^b(Q)$, called a Serre functor, such that
there are isomorphisms
\begin{equation} \label{Serreduality}
\Hom_{\CD^b(Q)}(Y,\nu X)\simeq D\Hom_{\CD^b(Q)}(X,Y)\simeq\Hom_{\CD^b(Q)}(\nu^-Y,X),
\end{equation}
which are natural in $X$ and $Y$, where $D$ is the duality $\Hom_\Bbbk(?,\Bbbk)$ on $\mod\Bbbk$ (see~\cite{RV}*{I}).
The Serre functor $\nu$ restricts to an equivalence between the subcategory $\CP_Q$ of the projective representations
in $\rep(Q)$ and the subcategory $\CI_Q$ of the injective representations in $\rep(Q)$.
This restriction is called a Nakayama functor.
For each vertex $a$ of $Q$, we denote by $P_a$ and $I_a$ the indecomposable projective
and injective representation in $\rep(Q)$ at $a$, respectively.
We note that up to isomorphism, these are the only indecomposable objects of $\CP_Q$ and $\CI_Q$, respectively.

The existence of a Serre functor $\nu$ is closely related to the existence of almost split triangles in $\CD^b(Q)$.
Namely, we consider the auto-equivalence $\tau=\nu\circ[-1]\simeq[-1]\circ\nu$ of $\CD^b(Q)$,
and call it the Auslander-Reiten translation.
Then there is an almost split triangle of the form $\tau C\to B\to C\to(\tau C)[1]$ for any indecomposable
object $C$ in $\CD^b(Q)$, and there is an almost split triangle of the form $A\to B'\to\tau^-A\to A[1]$
for any indecomposable object $A$ in $\CD^b(Q)$.

One defines the Grothendieck group $K_0(\rep(Q))$ of $\rep(Q)$ as
the quotients of the free abelian group with basis formed by the isomorphism classes $[X]$ of objects $X$ in $\rep(Q)$,
modulo the subgroup generated by $[A]-[B]+[C]$ for all short exact sequences $0\to A\to B\to C\to 0$ in $\rep(Q)$.
The group $K_0(\rep(Q))$ is isomorphic with $\BZ^{Q_0}$ via the map sending the class of $[X]$ to $\dimv X$.
We will treat this isomorphism as an identification.

The Grothendieck group $K_0(\CD^b(Q))$ of the category $\CD^b(Q)$ is defined in a similar way,
the only difference is that one replaces the sequences $0\to A\to B\to C\to 0$ by the triangles $A\to B\to C\to A[1]$
when forming the quotient.
The embedding of $\rep(Q)$ in $\CD^b(Q)$ induces a group isomorphism from $K_0(\rep(Q))$ to $K_0(\CD^b(Q))$,
which we will treat as identification.
In particular we will use notation $\dimv X$ for $X\in\CD^b(Q)$.
Note that $\dimv X[i]=(-1)^i\cdot\dimv X$ for any object $X\in\CD^b(Q)$ and any integer $i$.

Given two objects $X$ and $Y$ in $\CD^b(Q)$ we denote by $[X,Y]$ the dimension of $\Hom_{\CD^b(Q)}(X,Y)$.
Moreover, we set $[X,Y]^i=[X[-i],Y]=[X,Y[i]]$ for any integer $i$.
If $X$ and $Y$ belong to $\rep(Q)$ then $[X,Y]^i=0$ provided $i\not\in\{0,1\}$.
We define the bilinear form $\frakb_Q\colon\BZ^{Q_0}\times\BZ^{Q_0}\to\BZ$ by the formula
\[
\frakb_Q(\dd,\ee)=\sum_{a\in Q_0}\dd_a\cdot\ee_a-\sum_{\alpha\in Q_1}\dd_{s\alpha}\cdot\ee_{t\alpha},
\]
for all $\dd,\ee\in\BZ^{Q_0}$. Then
\[
\frakb_Q(\dimv X,\dimv Y)=\sum_{i\in\BZ}(-1)^i\cdot[X,Y]^i,
\]
for all $X,Y\in\CD^b(Q)$ (see~\cite{Hap}*{III.1}). In particular, if $X,Y\in\rep(Q)$, then
\[
\frakb_Q(\dimv X,\dimv Y)=[X,Y]^0-[X,Y]^1=\dim_\Bbbk\Hom_Q(X,Y)-\dim_\Bbbk\Ext^1_Q(X,Y).
\]
Observe that the quadratic form associated with $\frakb_Q$ coincides with $\frakq_\Delta$,
where $\Delta$ is the underlying graph of $Q$ and $\frakq_\Delta$ is the quadratic form introduced in \eqref{qDelta}.
Recall that Yoneda lemma states that if $a \in Q_0$ and $M \in \rep (Q)$, then $[P_a,M]$ is the $a$-th coordinate of $\dimv M$.
This easily implies that $\frakb_Q(\dimv P_a,\dd)=\dd_a$, for each $\dd\in \BZ^{Q_0}$.

We collect below few facts concerning indecomposable objects in $\CD^b(Q)$ under the assumption that
$Q$ is a Dynkin quiver.

\begin{lem} \label{only1nonzero}
Let $Q$ be a Dynkin quiver, and $X$ and $Y$ be indecomposable objects in $\CD^b(Q)$.
Then:
\begin{enumerate}
\item[\textup{(1)}] $X\simeq\tau^{n_X}P_{a_X}$ for a unique pair $(n_X,a_X)\in\BZ\times Q_0$.
\item[\textup{(2)}] $[X,X]=1$ and $\frakq_\Delta(\dimv X)=1$.
\item[\textup{(3)}] $[X,Y]^i$ is non-zero for at most one integer $i$.
\item[\textup{(4)}] $[X,Y]\leq\hh^\Delta_{a_X}$.
\end{enumerate}
\end{lem}

\begin{proof}
The first three properties are well known (see for example~\cite{Hap}).

(4). We may assume that $[X,Y]>0$.
Applying the automorphism $\tau^{-n_X}$ we get
\[
[X,Y]=[\tau^{n_X}P_{a_X},Y]=[P_{a_X},\tau^{-n_X}Y]=\frakb_Q(\dimv P_{a_X},\dimv\tau^{-n_X}Y),
\]
where the last equality follows from~(3) and the assumption $[X,Y]>0$.
As we observed above $\frakb_Q(\dimv P_{a_X},\dimv\tau^{-n_X}Y)$ is the $a_X$-th coordinate of the vector $\dimv\tau^{-n_X}Y$.
By (2), $\frakq_\Delta(\dimv\tau^{-n_X}Y)=1$, i.e.\ $\dimv\tau^{-n_X}Y$ is a root, and the claim follows, since $\hh^\Delta$ is the maximal root.
\end{proof}

\subsection{Mesh categories for Dynkin graphs}
\label{subsect mesh}

Throughout this subsection $\Delta=(\Delta_0,\Delta_1)$ is a Dynkin graph.
We say that two elements $(p,a)$ and $(q,b)$ of the product $\BZ\times\Delta_0$ are equivalent
provided the integer $(q-p)+\dist(a,b)$ is even.
This is an equivalence relation since $\Delta$ is a tree, and thus $\BZ\times\Delta_0$
is partitioned into two parts
\[
\BZ\times\Delta_0=(\BZ\times\Delta_0)^\vv\;\sqcup\;(\BZ\times\Delta_0)^\mm.
\]
In order to decide which part stands for $(\BZ\times\Delta_0)^\vv$, we choose
a base vertex $b_0\in\Delta_0$ and require that $(0,b_0)$ belongs to $(\BZ\times\Delta_0)^\vv$.

We define an infinite quiver $\BZ\Delta$ without multiple arrows as follows.
The set $(\BZ\Delta)_0$ of vertices of $\BZ \Delta$ consists of $\vv_{p,a}$, where $(p,a)\in(\BZ\times\Delta_0)^\vv$.
There is an arrow in $\BZ\Delta$ starting at $\vv_{p,a}$ and terminating at $\vv_{q,b}$ if and only if
$a$ and $b$ are adjacent in $\Delta$ and $q-p=1$.
For example, if $\Delta=\BA_4$:
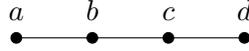
\[
\begin{tikzpicture}
\draw (0,0) -- (3,0);
\filldraw (0,0) circle (2pt) node [above=3pt] {$a$};
\filldraw (1,0) circle (2pt) node [above=3pt] {$b$};
\filldraw (2,0) circle (2pt) node [above=3pt] {$c$};
\filldraw (3,0) circle (2pt) node [above=3pt] {$d$};
\end{tikzpicture}
\]
and we choose $b_0=b$, then $\BZ\Delta$ has the form
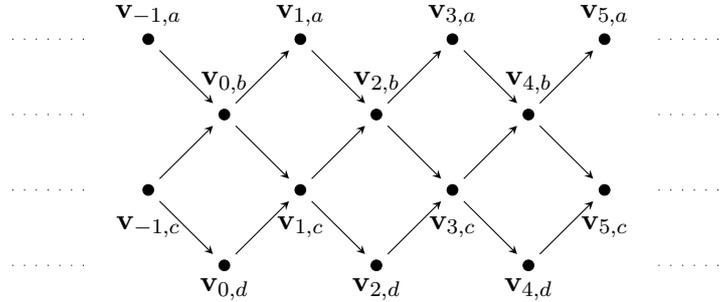
\[
\begin{tikzpicture} [>=stealth,shorten <= 6pt,shorten >=6pt] 
\filldraw
 (-1,3) circle (2pt) node [above=1pt] {$\vv_{-1,a}$}
 (1,3) circle (2pt) node [above=1pt] {$\vv_{1,a}$}
 (3,3) circle (2pt) node [above=1pt] {$\vv_{3,a}$}
 (5,3) circle (2pt) node [above=1pt] {$\vv_{5,a}$}
 (0,2) circle (2pt) node [above=4pt] {$\vv_{0,b}$}
 (2,2) circle (2pt) node [above=4pt] {$\vv_{2,b}$}
 (4,2) circle (2pt) node [above=4pt] {$\vv_{4,b}$}
 (-1,1) circle (2pt) node [below=6pt] {$\vv_{-1,c}$}
 (1,1) circle (2pt) node [below=6pt] {$\vv_{1,c}$}
 (3,1) circle (2pt) node [below=6pt] {$\vv_{3,c}$}
 (5,1) circle (2pt) node [below=6pt] {$\vv_{5,c}$}
 (0,0) circle (2pt) node [below=1pt] {$\vv_{0,d}$}
 (2,0) circle (2pt) node [below=1pt] {$\vv_{2,d}$}
 (4,0) circle (2pt) node [below=1pt] {$\vv_{4,d}$};
\draw [->] (-1,3) -- (0,2);
\draw [->] (0,2) -- (1,3);
\draw [->] (1,3) -- (2,2);
\draw [->] (2,2) -- (3,3);
\draw [->] (3,3) -- (4,2);
\draw [->] (4,2) -- (5,3);
\draw [->] (-1,1) -- (0,2);
\draw [->] (0,2) -- (1,1);
\draw [->] (1,1) -- (2,2);
\draw [->] (2,2) -- (3,1);
\draw [->] (3,1) -- (4,2);
\draw [->] (4,2) -- (5,1);
\draw [->] (-1,1) -- (0,0);
\draw [->] (0,0) -- (1,1);
\draw [->] (1,1) -- (2,0);
\draw [->] (2,0) -- (3,1);
\draw [->] (3,1) -- (4,0);
\draw [->] (4,0) -- (5,1);
\draw [loosely dotted] (-3,3) -- (-1.5,3);
\draw [loosely dotted] (-3,2) -- (-1.5,2);
\draw [loosely dotted] (-3,1) -- (-1.5,1);
\draw [loosely dotted] (-3,0) -- (-1.5,0);
\draw [loosely dotted] (5.5,3) -- (7,3);
\draw [loosely dotted] (5.5,2) -- (7,2);
\draw [loosely dotted] (5.5,1) -- (7,1);
\draw [loosely dotted] (5.5,0) -- (7,0);
\end{tikzpicture}
\]

For each pair $(p,a)$ in $(\BZ\times\Delta_0)^\mm$,
we consider the smallest subquiver $\mm_{p,a}$ of $\BZ\Delta$, called a mesh, containing all paths (of length two)
starting at $\vv_{p-1,a}$ and terminating at $\vv_{p+1,a}$:
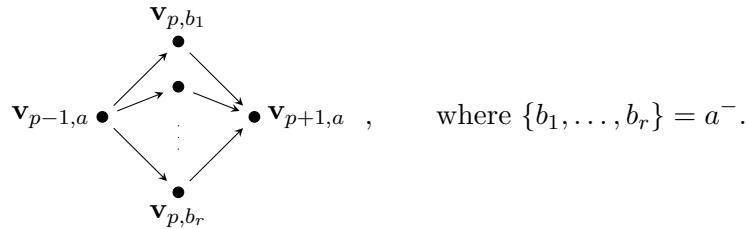
\[
\vcenter{\hbox{
\begin{tikzpicture} [>=stealth,shorten <= 6pt,shorten >=6pt]
\filldraw
 (-1,0) circle (2pt) node [left=1pt] {$\vv_{p-1,a}$}
 (0,1) circle (2pt) node [above=1pt] {$\vv_{p,b_1}$}
 (0,.4) circle (2pt)
 (0,-1) circle (2pt) node [below=1pt] {$\vv_{p,b_r}$}
 (1,0) circle (2pt) node [right=1pt] {$\vv_{p+1,a}$};
\draw [->] (-1,0) -- (0,1);
\draw [->] (0,1) -- (1,0);
\draw [->] (-1,0) -- (0,.4);
\draw [->] (0,.4) -- (1,0);
\draw [->] (-1,0) -- (0,-1);
\draw [->] (0,-1) -- (1,0);
\draw [loosely dotted] (0,.1) -- (0,-.7);
\end{tikzpicture}
}},\qquad
\text{where $\{b_1,\ldots,b_r\}=a^-$.}
\]
We denote by $(\BZ\Delta)_2$ the set of meshes $\mm_{p,a}$, where $(p,a)\in(\BZ\times\Delta_0)^\mm$.
With each mesh $\mm_{p,a}$ we associate its mesh relation, i.e.\ the sum of the paths
starting at $\vv_{p-1,a}$ and terminating at $\vv_{p+1,a}$, considered as morphisms
in the path category $\Bbbk[\BZ\Delta]$ of $\BZ\Delta$.
The mesh category $\Bbbk(\BZ\Delta)$ of $\BZ\Delta$ is the quotient of $\Bbbk[\BZ\Delta]$
modulo the ideal generated by all mesh relations.

Now let $Q$ be a Dynkin quiver with underlying graph $\Delta$.
An important fact is that $\Bbbk(\BZ\Delta)$ is equivalent as a $\Bbbk$-linear category
to the category of indecomposable objects in $\CD^b(Q)$.
When $Q$ is fixed, then we shall identify $(\BZ\Delta)_0$ with a complete set of pairwise
non-isomorphic indecomposable objects of $\CD^b(Q)$.
Moreover, we may also assume that under this identification $P_a=\vv_{p_a,a}$ with an appropriate
integer $p_a$, for each vertex $a \in Q_0 = \Delta_0$.
The three crucial auto-equivalences of $\CD^b(Q)$: the Auslander-Reiten translation $\tau$,
the Serre functor $\nu$ and the shift functor $[1]$ act on the indecomposable objects by the formulas
\[
\tau(\vv_{p,a})=\vv_{p-2,a},
\qquad
\nu(\vv_{p,a})=\vv_{p+n_\Delta-2,\phi_\Delta(a)},
\qquad
\vv_{p,a}[1]=\vv_{p+n_\Delta,\phi_\Delta(a)},
\]
where $n_\Delta$ was defined in~\eqref{nDelta} and $\phi_\Delta$ is the automorphism of $\Delta$ defined as follows:
$\phi_\Delta$ is the unique non-trivial involution of $\Delta$ provided
$\Delta$ is either of type $\BA_n$ with $n\geq 2$, or $\BD_n$ with $n$ odd, or $\BE_6$;
and $\phi_\Delta$ is the identity on $\Delta$ for the remaining Dynkin graphs.
We note that the automorphism of $\BZ\Delta$ induced by $\nu$ is sometimes called a Nakayama permutation (\cite{Gab}*{6.5}).
We also remark that the quiver $\BZ\Delta$ is isomorphic to the quiver $\BZ Q$ defined in~\cite{Rdarst}.

The almost split triangles in $\CD^b(Q)$ are parameterized by the meshes in $\BZ\Delta$.
More precisely, there is an almost split triangle of the form
\[
\AR(\mm_{p,a})\colon\vv_{p-1,a}\to\bigoplus_{b\in a^-}\vv_{p,b}\to\vv_{p+1,a}\to\vv_{p-1,a}[1],
\]
where $\vv_{p-1,a}[1]=\vv_{p+n_\Delta-1,\phi_\Delta(a)}$,
for any mesh $\mm_{p,a}\in(\BZ\Delta)_2$.

We note that the paths in $\BZ\Delta$ have the form
\[
\omega\colon\quad\vv_{p,a_p}\to\vv_{p+1,a_{p+1}}\to\quad\cdots\quad\to\vv_{q-1,a_{q-1}}\to\vv_{q,a_q},
\]
where each two consecutive vertices in the sequence $(a_p,a_{p+1},\ldots,a_{q-1},a_q)$ are adjacent.
Therefore $\omega$ can be viewed as a lifting of a walk in $\Delta$.
In particular, $q-p\geq\dist(a_p,a_q)$.
The path $\omega$ is called sectional if $a_{i-1}\neq a_{i+1}$ for any integer $i$ with $p<i<q$.
Since $\Delta$ is a tree, this condition is equivalent to the fact that the vertices $a_p,\ldots,a_q$ are pairwise different,
and also equivalent to the equality $q-p=\dist(a_p,a_q)$.

\begin{lem} \label{homifpath}
Let $\vv_{p,a}$ and $\vv_{q,b}$ be vertices in $\BZ\Delta$.
Then:
\begin{enumerate}
\item[\textup{(1)}] $[\vv_{p,a},\vv_{p,a}]=1$.
\item[\textup{(2)}] $[\vv_{p,a},\vv_{q,b}]=[\vv_{q,b},\vv_{p+n_\Delta-2,\phi_\Delta(a)}]$.
\item[\textup{(3)}] If $[\vv_{p,a},\vv_{q,b}]>0$ then $p+\dist(a,b)\leq q\leq p+n_\Delta-2-\dist(b,\phi_\Delta(a))$.
\item[\textup{(4)}] $[\vv_{p,a},\vv_{q,b}]\leq\min(\hh^\Delta_a,\hh^\Delta_b)$
\end{enumerate}
\end{lem}

\begin{proof}
(1) follows from Lemma~\ref{only1nonzero}(2), but can also be derived directly from the definition
of the mesh category $\Bbbk(\BZ\Delta)$.
(2) is a consequence of the Serre duality \eqref{Serreduality}.

(3). If $[\vv_{p,a},\vv_{q,b}]>0$ then also $[\vv_{q,b},\vv_{p+n_\Delta-2,\phi_\Delta(a)}]>0$, by (2).
Hence there are paths in $\BZ\Delta$ from $\vv_{p,a}$ to $\vv_{q,b}$ and from $\vv_{q,b}$ to $\vv_{p+n_\Delta-2,\phi_\Delta(a)}$.

(4). Since $P_a=\vv_{p_a,a}$, $\vv_{p,a}=\tau^rP_a$ for some integer $r$.
By Lemma~\ref{only1nonzero}(4), $[\vv_{p,a},\vv_{q,b}]\leq\hh^\Delta_a$.
The other inequality $[\vv_{p,a},\vv_{q,b}]\leq\hh^\Delta_b$ follows from the first one and (2).
\end{proof}

Let us explain how using the above lemma and almost split sequences, we can calculate the dimension
$[\vv_{p,a},\vv_{q,b}]$ for all vertices $\vv_{p,a}$ and $\vv_{q,b}$,
Namely, if $q\leq p$ or $q\geq p+n_\Delta-1$ then $[\vv_{p,a},\vv_{q,b}]=0$ except $[\vv_{p,a},\vv_{p,a}]=1$.
We obtain formulas in the remaining cases by induction on $q$, using the following lemma.

\begin{lem} \label{homalgorithm}
Let $\vv_{p,a}$ and $\vv_{q,b}$ be vertices in $\BZ\Delta$ such that $p<q<p+n_\Delta$.
Then
\[
[\vv_{p,a},\vv_{q,b}]=\sum_{c\in b^-}[\vv_{p,a},\vv_{q-1,c}]-[\vv_{p,a},\vv_{q-2,b}].
\]
\end{lem}

\begin{proof}
Applying the functor $\Hom_{\CD^b(Q)}(\vv_{p,a},?)$
to the triangle $\AR(\mm_{q-1,b})$ we get the exact sequence
\begin{multline*}
\ldots\to\Hom(\vv_{p,a},\vv_{q-n_\Delta,\phi_\Delta(b)})\to\Hom(\vv_{p,a},\vv_{q-2,b})\to\bigoplus_{c\in b^-}\Hom(\vv_{p,a},\vv_{q-1,c})\to \\
\to\Hom(\vv_{p,a},\vv_{q,b})\to\Hom(\vv_{p,a},\vv_{q+n_\Delta-2,\phi_\Delta(b)})\to\ldots
\end{multline*}
The two extreme homomorphism spaces are zero by Lemma~\ref{homifpath}(3), and the claim follows.
\end{proof}

Applying the above to sectional paths we get the following.

\begin{lem} \label{hom1section}
If $\vv_{p,a_p}\to\vv_{p+1,a_{p+1}}\to\ldots\to\vv_{q-1,a_{q-1}}\to\vv_{q,a_q}$
is a sectional path in $\BZ\Delta$ then $[\vv_{p,a_p},\vv_{q,a_q}]=1$.
\end{lem}

\begin{proof}
The claim follows by induction on the length $(q-p)$ of the path,
where the base step $q-p=0$ follows from Lemma~\ref{homifpath}.
For the induction step we apply Lemma~\ref{homalgorithm} for $\vv_{p,a}=\vv_{p,a_p}$ and $\vv_{q,b}=\vv_{q,a_q}$,
and use that there is no path in $\BZ\Delta$ from $\vv_{p,a_p}$ to $\vv_{q-2,a_q}$, and if there is a path from
$\vv_{p,a_p}$ to $\vv_{q-1,c}$ with $c\in(a_q)^-$ then $c=a_{q-1}$.
\end{proof}

\subsection{Defect functions on meshes}

Throughout this subsection $Q$ is a Dynkin quiver with the underlying graph $\Delta$.
In paricular, $Q_0=\Delta_0$.
Consider a triangle
\[
\sigma\colon A\xrightarrow{\alpha}B\xrightarrow{\beta}C\xrightarrow{\gamma}A[1]
\]
and an object $X$ in $\CD^b(Q)$.
There is a commutative diagram of the form
\[
\xymatrixcolsep{7pc}
\xymatrix{
\Hom_{\CD^b(Q)}(X,C)\ar[r]^-{\Hom_{\CD^b(Q)}(X,\gamma)}\ar[d]&\Hom_{\CD^b(Q)}(X,A[1])\ar[d]\\
D\Hom_{\CD^b(Q)}(C[-1],\tau X)\ar[r]^-{D\Hom_{\CD^b(Q)}(\gamma[-1],\tau X)}&D\Hom_{\CD^b(Q)}(A,\tau X)
}
\]
where the vertical arrows represent $\Bbbk$-linear isomorphisms obtained by applying
the Auslander-Reiten translation $\tau$ and the Serre duality \eqref{Serreduality}.
This inspires to define the following integer-valued function measuring how far is a triangle from being split.

\begin{df} \label{defdeltasigma}
Given a triangle $\sigma\colon A\xrightarrow{\alpha}B\xrightarrow{\beta}C\xrightarrow{\gamma}A[1]$ in $\CD^b(Q)$
we define a non-negative function
\[
\delta_\sigma\colon(\BZ\Delta)_2\to\BZ,
\qquad
\delta_\sigma(\mm_{p,a})=\rk\bigl(\Hom_{\CD^b(Q)}(\vv_{p+1,a},\gamma)\bigr)=\rk\bigl(\Hom_{\CD^b(Q)}(\gamma[-1],\vv_{p-1,a})\bigr).
\]
\end{df}

Some fundamental properties of $\delta_\sigma$ being easy consequences of
the definition are collected in the following lemma.

\begin{lem} \label{fundamental}
Let $\sigma\colon A\xrightarrow{\alpha}B\xrightarrow{\beta}C\xrightarrow{\gamma}A[1]$ be a triangle in $\CD^b(Q)$.
Then the following hold:
\begin{enumerate}
\item[\textup{(1)}] $\sigma$ splits if and only if $\delta_\sigma=0$.
\item[\textup{(2)}] If $B=0$, then $\delta_\sigma(\mm_{p,a})=[A,\vv_{p-1,a}]=[\vv_{p+1,a},C]$.
\item[\textup{(3)}] $\delta_\sigma(\mm_{p,a})\leq[A,\vv_{p-1,a}]$ and $\delta_\sigma(\mm_{p,a})\leq[\vv_{p+1,a},C]$.
\item[\textup{(4)}] $\delta_{\sigma'}\leq\delta_\sigma$ if $\sigma'$ is a pullback or a pushout of $\sigma$.
\item[\textup{(5)}] If $\sigma=\AR(\mm_{q,b})$, then
\[
\delta_\sigma(\mm_{p,a})=\begin{cases}
1,&(p,a)=(q,b),\\
0,&\text{otherwise}.
\end{cases} \eqno \qed
\]
\end{enumerate}
\end{lem}

\begin{ex}
Let $\Delta$ be the Dynkin graph $\BD_6$:
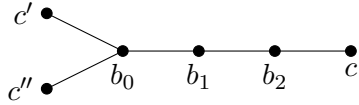
\[
\begin{tikzpicture}
\draw (0,.5) -- (1,0);
\draw (0,-.5) -- (1,0);
\draw (1,0) -- (4,0);
\filldraw
 (0,.5) circle (2pt) node [left=1pt] {$c'$}
 (0,-.5) circle (2pt) node [left=1pt] {$c''$}
 (1,0) circle (2pt) node [below=1pt] {$b_0$}
 (2,0) circle (2pt) node [below=1pt] {$b_1$}
 (3,0) circle (2pt) node [below=1pt] {$b_2$}
 (4,0) circle (2pt) node [below=1pt] {$c$};
\end{tikzpicture}
\]
Thus $n_\Delta=10$ and the quiver $\BZ\Delta$ has the form
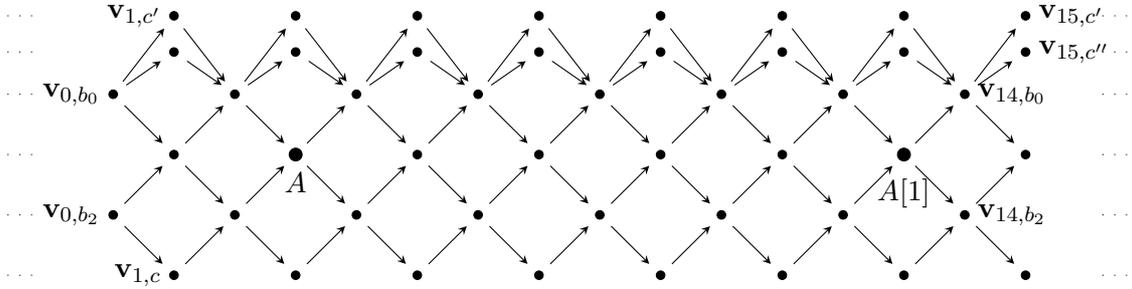
\[
\begin{tikzpicture} [>=stealth,shorten <= 6pt,shorten >=6pt,scale=.8] 
\filldraw
 (1,4.3) circle (2pt) node [left=1pt] {$\vv_{1,c'}$}
 (3,4.3) circle (2pt)
 (5,4.3) circle (2pt)
 (7,4.3) circle (2pt)
 (9,4.3) circle (2pt)
 (11,4.3) circle (2pt)
 (13,4.3) circle (2pt)
 (15,4.3) circle (2pt) node [right=1pt] {$\vv_{15,c'}$}
  (1,3.7) circle (2pt)
 (3,3.7) circle (2pt)
 (5,3.7) circle (2pt)
 (7,3.7) circle (2pt)
 (9,3.7) circle (2pt)
 (11,3.7) circle (2pt)
 (13,3.7) circle (2pt)
 (15,3.7) circle (2pt) node [right=1pt] {$\vv_{15,c''}$}
 (0,3) circle (2pt) node [left=1pt] {$\vv_{0,b_0}$}
 (2,3) circle (2pt)
 (4,3) circle (2pt)
 (6,3) circle (2pt)
 (8,3) circle (2pt)
 (10,3) circle (2pt)
 (12,3) circle (2pt)
 (14,3) circle (2pt) node [right=1pt] {$\vv_{14,b_0}$}
 (1,2) circle (2pt)
 (3,2) circle (3pt) node [below=3pt] {$A$}
 (5,2) circle (2pt)
 (7,2) circle (2pt)
 (9,2) circle (2pt)
 (11,2) circle (2pt)
 (13,2) circle (3pt) node [below=5pt] {$A[1]$}
 (15,2) circle (2pt)
 (0,1) circle (2pt) node [left=1pt] {$\vv_{0,b_2}$}
 (2,1) circle (2pt)
 (4,1) circle (2pt)
 (6,1) circle (2pt)
 (8,1) circle (2pt)
 (10,1) circle (2pt)
 (12,1) circle (2pt)
 (14,1) circle (2pt) node [right=1pt] {$\vv_{14,b_2}$}
 (1,0) circle (2pt) node [left=1pt] {$\vv_{1,c}$}
 (3,0) circle (2pt)
 (5,0) circle (2pt)
 (7,0) circle (2pt)
 (9,0) circle (2pt)
 (11,0) circle (2pt)
 (13,0) circle (2pt)
 (15,0) circle (2pt);
\draw [->] (1,4.3) -- (2,3);
\draw [->] (3,4.3) -- (4,3);
\draw [->] (5,4.3) -- (6,3);
\draw [->] (7,4.3) -- (8,3);
\draw [->] (9,4.3) -- (10,3);
\draw [->] (11,4.3) -- (12,3);
\draw [->] (13,4.3) -- (14,3);
\draw [->] (0,3) -- (1,4.3);
\draw [->] (2,3) -- (3,4.3);
\draw [->] (4,3) -- (5,4.3);
\draw [->] (6,3) -- (7,4.3);
\draw [->] (8,3) -- (9,4.3);
\draw [->] (10,3) -- (11,4.3);
\draw [->] (12,3) -- (13,4.3);
\draw [->] (14,3) -- (15,4.3);
\draw [->] (1,3.7) -- (2,3);
\draw [->] (3,3.7) -- (4,3);
\draw [->] (5,3.7) -- (6,3);
\draw [->] (7,3.7) -- (8,3);
\draw [->] (9,3.7) -- (10,3);
\draw [->] (11,3.7) -- (12,3);
\draw [->] (13,3.7) -- (14,3);
\draw [->] (0,3) -- (1,3.7);
\draw [->] (2,3) -- (3,3.7);
\draw [->] (4,3) -- (5,3.7);
\draw [->] (6,3) -- (7,3.7);
\draw [->] (8,3) -- (9,3.7);
\draw [->] (10,3) -- (11,3.7);
\draw [->] (12,3) -- (13,3.7);
\draw [->] (14,3) -- (15,3.7);
\draw [->] (1,2) -- (2,3);
\draw [->] (3,2) -- (4,3);
\draw [->] (5,2) -- (6,3);
\draw [->] (7,2) -- (8,3);
\draw [->] (9,2) -- (10,3);
\draw [->] (11,2) -- (12,3);
\draw [->] (13,2) -- (14,3);
\draw [->] (0,3) -- (1,2);
\draw [->] (2,3) -- (3,2);
\draw [->] (4,3) -- (5,2);
\draw [->] (6,3) -- (7,2);
\draw [->] (8,3) -- (9,2);
\draw [->] (10,3) -- (11,2);
\draw [->] (12,3) -- (13,2);
\draw [->] (14,3) -- (15,2);
\draw [->] (1,2) -- (2,1);
\draw [->] (3,2) -- (4,1);
\draw [->] (5,2) -- (6,1);
\draw [->] (7,2) -- (8,1);
\draw [->] (9,2) -- (10,1);
\draw [->] (11,2) -- (12,1);
\draw [->] (13,2) -- (14,1);
\draw [->] (0,1) -- (1,2);
\draw [->] (2,1) -- (3,2);
\draw [->] (4,1) -- (5,2);
\draw [->] (6,1) -- (7,2);
\draw [->] (8,1) -- (9,2);
\draw [->] (10,1) -- (11,2);
\draw [->] (12,1) -- (13,2);
\draw [->] (14,1) -- (15,2);
\draw [->] (1,0) -- (2,1);
\draw [->] (3,0) -- (4,1);
\draw [->] (5,0) -- (6,1);
\draw [->] (7,0) -- (8,1);
\draw [->] (9,0) -- (10,1);
\draw [->] (11,0) -- (12,1);
\draw [->] (13,0) -- (14,1);
\draw [->] (0,1) -- (1,0);
\draw [->] (2,1) -- (3,0);
\draw [->] (4,1) -- (5,0);
\draw [->] (6,1) -- (7,0);
\draw [->] (8,1) -- (9,0);
\draw [->] (10,1) -- (11,0);
\draw [->] (12,1) -- (13,0);
\draw [->] (14,1) -- (15,0);
\draw [loosely dotted] (-2,4.3) -- (-1,4.3);
\draw [loosely dotted] (-2,3.7) -- (-1,3.7);
\draw [loosely dotted] (-2,3) -- (-1,3);
\draw [loosely dotted] (-2,2) -- (-1,2);
\draw [loosely dotted] (-2,1) -- (-1,1);
\draw [loosely dotted] (-2,0) -- (-1,0);
\draw [loosely dotted] (16,4.3) -- (17,4.3);
\draw [loosely dotted] (16,3.7) -- (17,3.7);
\draw [loosely dotted] (16,3) -- (17,3);
\draw [loosely dotted] (16,2) -- (17,2);
\draw [loosely dotted] (16,1) -- (17,1);
\draw [loosely dotted] (16,0) -- (17,0);
\end{tikzpicture}
\]
We consider the triangle $\sigma\colon A\to 0\to A[1]\xrightarrow{1}A[1]$ with $A=\vv_{3,b_1}$.
Then $A[1]=\vv_{13,b_1}$ and $\delta_\sigma(\mm_{q,b})=[\vv_{3,b_1},\vv_{q-1,b}]=[\vv_{q+1,b},\vv_{13,b_1}]$, by Lemma~\ref{fundamental}(2).
We find $\delta_\sigma$ by calculating the dimensions $[\vv_{3,b_1},\vv_{q-1,b}]$, which can be done by the method based on Lemma~\ref{homalgorithm}.
We illustrate the function $\delta_\sigma$ by writing each non-zero value
$\delta_\sigma(\mm_{p,a})$ between the vertices $\vv_{p-1,a}$ and $\vv_{p+1,a}$:
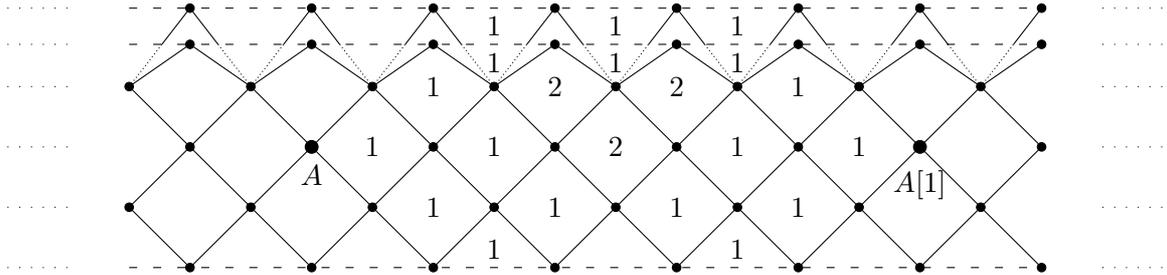
\[
\begin{tikzpicture} [>=stealth,scale=.8]
\filldraw
 (1,4.3) circle (2pt)
 (3,4.3) circle (2pt)
 (5,4.3) circle (2pt)
 (7,4.3) circle (2pt)
 (9,4.3) circle (2pt)
 (11,4.3) circle (2pt)
 (13,4.3) circle (2pt)
 (15,4.3) circle (2pt)
  (1,3.7) circle (2pt)
 (3,3.7) circle (2pt)
 (5,3.7) circle (2pt)
 (7,3.7) circle (2pt)
 (9,3.7) circle (2pt)
 (11,3.7) circle (2pt)
 (13,3.7) circle (2pt)
 (15,3.7) circle (2pt)
 (0,3) circle (2pt)
 (2,3) circle (2pt)
 (4,3) circle (2pt)
 (6,3) circle (2pt)
 (8,3) circle (2pt)
 (10,3) circle (2pt)
 (12,3) circle (2pt)
 (14,3) circle (2pt)
 (1,2) circle (2pt)
 (3,2) circle (3pt) node [below=3pt] {$A$}
 (5,2) circle (2pt)
 (7,2) circle (2pt)
 (9,2) circle (2pt)
 (11,2) circle (2pt)
 (13,2) circle (3pt) node [below=5pt] {$A[1]$}
 (15,2) circle (2pt)
 (0,1) circle (2pt)
 (2,1) circle (2pt)
 (4,1) circle (2pt)
 (6,1) circle (2pt)
 (8,1) circle (2pt)
 (10,1) circle (2pt)
 (12,1) circle (2pt)
 (14,1) circle (2pt)
 (1,0) circle (2pt)
 (3,0) circle (2pt)
 (5,0) circle (2pt)
 (7,0) circle (2pt)
 (9,0) circle (2pt)
 (11,0) circle (2pt)
 (13,0) circle (2pt)
 (15,0) circle (2pt);
\draw
 node at (6,4) {$1$}
 node at (8,4) {$1$}
 node at (10,4) {$1$}
 node at (6,3.4) {$1$}
 node at (8,3.4) {$1$}
 node at (10,3.4) {$1$}
 node at (5,3) {$1$}
 node at (7,3) {$2$}
 node at (9,3) {$2$}
 node at (11,3) {$1$}
 node at (4,2) {$1$}
 node at (6,2) {$1$}
 node at (8,2) {$2$}
 node at (10,2) {$1$}
 node at (12,2) {$1$}
 node at (5,1) {$1$}
 node at (7,1) {$1$}
 node at (9,1) {$1$}
 node at (11,1) {$1$}
 node at (6,.3) {$1$}
 node at (10,.3) {$1$};
\draw [loosely dashed] (0,4.3) -- (15,4.3);
\draw [loosely dashed] (0,3.7) -- (15,3.7);
\draw [loosely dashed] (0,0) -- (15,0);
\draw (1,4.3) -- (1+6/13,3.7);
\draw [densely dotted] (1+6/13,3.7) -- (2,3);
\draw (3,4.3) -- (3+6/13,3.7);
\draw [densely dotted] (3+6/13,3.7) -- (4,3);
\draw (5,4.3) -- (5+6/13,3.7);
\draw [densely dotted] (5+6/13,3.7) -- (6,3);
\draw (7,4.3) -- (7+6/13,3.7);
\draw [densely dotted] (7+6/13,3.7) -- (8,3);
\draw (9,4.3) -- (9+6/13,3.7);
\draw [densely dotted] (9+6/13,3.7) -- (10,3);
\draw (11,4.3) -- (11+6/13,3.7);
\draw [densely dotted] (11+6/13,3.7) -- (12,3);
\draw (13,4.3) -- (13+6/13,3.7);
\draw [densely dotted] (13+6/13,3.7) -- (14,3);
\draw [densely dotted] (0,3) -- (0+7/13,3.7);
\draw (0+7/13,3.7) -- (1,4.3);
\draw [densely dotted] (2,3) -- (2+7/13,3.7);
\draw (2+7/13,3.7) -- (3,4.3);
\draw [densely dotted] (4,3) -- (4+7/13,3.7);
\draw (4+7/13,3.7) -- (5,4.3);
\draw [densely dotted] (6,3) -- (6+7/13,3.7);
\draw (6+7/13,3.7) -- (7,4.3);
\draw [densely dotted] (8,3) -- (8+7/13,3.7);
\draw (8+7/13,3.7) -- (9,4.3);
\draw [densely dotted] (10,3) -- (10+7/13,3.7);
\draw (10+7/13,3.7) -- (11,4.3);
\draw [densely dotted] (12,3) -- (12+7/13,3.7);
\draw (12+7/13,3.7) -- (13,4.3);
\draw [densely dotted] (14,3) -- (14+7/13,3.7);
\draw (14+7/13,3.7) -- (15,4.3);
\draw  (1,3.7) -- (2,3);
\draw  (3,3.7) -- (4,3);
\draw  (5,3.7) -- (6,3);
\draw  (7,3.7) -- (8,3);
\draw  (9,3.7) -- (10,3);
\draw  (11,3.7) -- (12,3);
\draw  (13,3.7) -- (14,3);
\draw  (0,3) -- (1,3.7);
\draw  (2,3) -- (3,3.7);
\draw  (4,3) -- (5,3.7);
\draw  (6,3) -- (7,3.7);
\draw  (8,3) -- (9,3.7);
\draw  (10,3) -- (11,3.7);
\draw  (12,3) -- (13,3.7);
\draw  (14,3) -- (15,3.7);
\draw  (1,2) -- (2,3);
\draw  (3,2) -- (4,3);
\draw  (5,2) -- (6,3);
\draw  (7,2) -- (8,3);
\draw  (9,2) -- (10,3);
\draw  (11,2) -- (12,3);
\draw  (13,2) -- (14,3);
\draw  (0,3) -- (1,2);
\draw  (2,3) -- (3,2);
\draw  (4,3) -- (5,2);
\draw  (6,3) -- (7,2);
\draw  (8,3) -- (9,2);
\draw  (10,3) -- (11,2);
\draw  (12,3) -- (13,2);
\draw  (14,3) -- (15,2);
\draw  (1,2) -- (2,1);
\draw  (3,2) -- (4,1);
\draw  (5,2) -- (6,1);
\draw  (7,2) -- (8,1);
\draw  (9,2) -- (10,1);
\draw  (11,2) -- (12,1);
\draw  (13,2) -- (14,1);
\draw  (0,1) -- (1,2);
\draw  (2,1) -- (3,2);
\draw  (4,1) -- (5,2);
\draw  (6,1) -- (7,2);
\draw  (8,1) -- (9,2);
\draw  (10,1) -- (11,2);
\draw  (12,1) -- (13,2);
\draw  (14,1) -- (15,2);
\draw  (1,0) -- (2,1);
\draw  (3,0) -- (4,1);
\draw  (5,0) -- (6,1);
\draw  (7,0) -- (8,1);
\draw  (9,0) -- (10,1);
\draw  (11,0) -- (12,1);
\draw  (13,0) -- (14,1);
\draw  (0,1) -- (1,0);
\draw  (2,1) -- (3,0);
\draw  (4,1) -- (5,0);
\draw  (6,1) -- (7,0);
\draw  (8,1) -- (9,0);
\draw  (10,1) -- (11,0);
\draw  (12,1) -- (13,0);
\draw  (14,1) -- (15,0);
\draw [loosely dotted] (-2,4.3) -- (-1,4.3);
\draw [loosely dotted] (-2,3.7) -- (-1,3.7);
\draw [loosely dotted] (-2,3) -- (-1,3);
\draw [loosely dotted] (-2,2) -- (-1,2);
\draw [loosely dotted] (-2,1) -- (-1,1);
\draw [loosely dotted] (-2,0) -- (-1,0);
\draw [loosely dotted] (16,4.3) -- (17,4.3);
\draw [loosely dotted] (16,3.7) -- (17,3.7);
\draw [loosely dotted] (16,3) -- (17,3);
\draw [loosely dotted] (16,2) -- (17,2);
\draw [loosely dotted] (16,1) -- (17,1);
\draw [loosely dotted] (16,0) -- (17,0);
\end{tikzpicture}
\]
Here we replace the arrows by edges, and additionally draw dashed segments between $\vv_{p-1,a}$ and $\vv_{p+1,a}$
if the vertex $a$ has degree $1$ in $\Delta$.
\end{ex}

We set $\langle X,Y\rangle=\sum_{i\leq 0}(-1)^i\cdot[X,Y]^i$ for any objects $X$ and $Y$ in $\CD^b(Q)$.
Our next aim is to define integer-valued functions on the set of meshes, using $\langle X,Y\rangle$.
We derive from Lemma~\ref{homifpath}(3) the following fact.

\begin{cor} \label{newhomifpath}
Assume $(p,a)$ and $(q,b)$ belong to $(\BZ\times\Delta_0)^\vv$.
If $\langle\vv_{p,a},\vv_{q,b}\rangle \neq 0$ then $q-p\geq\dist(a,b)$.
\qed
\end{cor}

\begin{lem}
Let $X,M,N\in\CD^b(Q)$ with $\dimv M=\dimv N$.
Then
\[
\langle X,N\rangle-\langle X,M\rangle=\langle N,\tau X\rangle-\langle M,\tau X\rangle.
\]
\end{lem}

\begin{proof}
By \eqref{Serreduality}, $[N,\tau X]^i=[N,\tau(X[i])]=[\nu^{-1}\tau(X[i]),N]=[X[i-1],N]=[X,N]^{1-i}$, for any integer $i$.
Consequently,
\[
\langle X,N\rangle-\langle N,\tau X\rangle=\frakb_Q(\dimv X,\dimv N)
=\frakb_Q(\dimv X,\dimv M)=\langle X,M\rangle-\langle M,\tau X\rangle.
\qedhere
\]
\end{proof}

\begin{df}
Let $M, N \in \CD^b(Q)$ be such that $\dimv M=\dimv N$.
We define the function $\delta_{M,N}\colon(\BZ\Delta)_2\to\BZ$ by
\[
\delta_{M,N}(\mm_{p,a})
=\langle\vv_{p+1,a},N\rangle-\langle\vv_{p+1,a},M\rangle
=\langle N,\vv_{p-1,a}\rangle-\langle M,\vv_{p-1,a}\rangle,
\]
for any mesh $\mm_{p,a}\in(\BZ\Delta)_2$.
\end{df}

By Corollary~\ref{newhomifpath} we conclude the following fact.

\begin{cor} \label{finitesupport}
Let $M,N\in\CD^b(Q)$ with $\dimv M=\dimv N$.
Then $\delta_{M,N}(\mm_{p,a})\neq 0$ only for finitely many
meshes $\mm_{p,a}\in(\BZ\Delta)_2$.
\qed
\end{cor}

Applying Hom functors we get the following.

\begin{cor} \label{deltasigmaalter}
$\delta_\sigma=\delta_{B,A\oplus C}$ for any triangle $\sigma\colon A\to B\to C\to A[1]$ in $\CD^b(Q)$.
\qed
\end{cor}

Combining the above corollary and Lemma~\ref{fundamental}(5) we get the following fact.

\begin{lem} \label{lemm213}
Let $(p,a)$ and $(q,b)$ belong to $(\BZ\times\Delta_0)^\vv$.
Then
\begin{multline*}
\langle\vv_{p,a}\oplus\vv_{p+2,a},\vv_{q,b}\rangle-\langle\bigoplus_{c\in a^-}\vv_{p+1,c},\vv_{q,b}\rangle
\\
=\langle\vv_{q,b},\vv_{p-2,a}\oplus\vv_{p,a}\rangle-\langle\vv_{q,b},\bigoplus_{c\in a^-}\vv_{p-1,c}\rangle
=\begin{cases}
1,&(q,b)=(p,a),\\
0,&\text{otherwise}.
\end{cases}
\tag*\qed
\end{multline*}

\end{lem}

Let $N$ and $X$ be objects of $\CD^b(Q)$ and assume that $X$ is indecomposable.
We denote by $\mult_X(N)$ the multiplicity of $X$ as a direct summand of $N$.
In particular,
\[
N\simeq\bigoplus_{(p,a)\in(\BZ\times\Delta_0)^\vv}(\vv_{p,a})^{\mult_{\vv_{p,a}}(N)}.
\]
As an immediate consequence of Lemma~\ref{lemm213} we get:

\begin{cor}
For any object $N$ of $\CD^b(Q)$ and $(p,a)\in(\BZ\times\Delta_0)^\vv$
\begin{align*}
\mult_{\vv_{p,a}}(N)&=\langle\vv_{p,a},N\rangle-\sum_{b\in a^-}\langle\vv_{p+1,b},N\rangle+\langle\vv_{p+2,a},N\rangle
\\
&=\langle N,\vv_{p-2,a}\rangle-\sum_{b\in a^-}\langle N,\vv_{p-1,b}\rangle+\langle N,\vv_{p,a}\rangle. \tag*\qed
\end{align*}
\end{cor}

\begin{cor} \label{multMN}
Let $M,N\in\CD^b(Q)$ with $\dimv M=\dimv N$ and $\vv_{p,a}\in(\BZ\Delta)_0$.
Then
\[
\mult_{\vv_{p,a}}(N)-\mult_{\vv_{p,a}}(M)=\delta_{M,N}(\mm_{p-1,a})-\sum_{b\in a^-}\delta_{M,N}(\mm_{p,b})+\delta_{M,N}(\mm_{p+1,a}).
\eqno \qed
\]
\end{cor}

Applying the above corollary for the vertices lying on a sectional path we obtain the following fact.

\begin{cor} \label{multMNsection}
Let $M,N\in\CD^b(Q)$ with $\dimv M=\dimv N$ and
\[
\vv_{p,a_p}\to\vv_{p+1,a_{p+1}}\to\quad\cdots\quad\to\vv_{q-1,a_{q-1}}\to\vv_{q,a_q}
\]
be a sectional path in $\BZ\Delta$.
Let $\CM$ be the subset of $(\BZ\times\Delta_0)^\mm$ consisting of pairs $(j,b)$ such that
$p\leq j\leq q$ and $b$ is adjacent to $a_j$, but does not belong to the set $\{a_p,\ldots,a_q\}$.
Then
\[
\sum_{i=p}^q\bigl(\mult_{\vv_{i,a_i}}(N)-\mult_{\vv_{i,a_i}}(M)\bigr)=\delta_{M,N}(\mm_{p-1,a_p})
-\sum_{(j,b)\in\CM}\delta_{M,N}(\mm_{j,b})
+\delta_{M,N}(\mm_{q+1,a_q}).
\eqno \qed
\]
\end{cor}

Observe in the above situation that if $b$ is adjacent to $a_j$, then $b$ belongs to the set $\{a_p,\ldots,a_q\}$
if and only if either $j>p$ and $b=a_{j-1}$ or $j<q$ and $b=a_{j+1}$.

Given a non-negative function $\delta\colon(\BZ\Delta)_2\to\BZ$, for instance $\delta_\sigma$ for a triangle $\sigma$,
we define its support
\[
\supp(\delta)=\{\mm\in(\BZ\Delta)_2;\;\delta(\mm)>0\}.
\]

\subsection{Application to type $\BD$}

Throughout this subsection $Q$ is a Dynkin quiver of type $\BD_n$, $n\geq 4$, with the underlying graph $\Delta$:
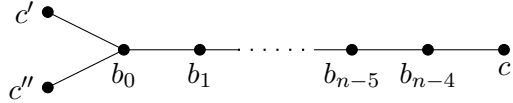
\[
\vcenter{\hbox{
\begin{tikzpicture}
\draw (0,.5) -- (1,0);
\draw (0,-.5) -- (1,0);
\draw (1,0) -- (2.5,0);
\draw [loosely dotted,thick] (2.5,0) -- (3.5,0);
\draw (3.5,0) -- (6,0);
\filldraw
 (0,.5) circle (2pt) node [left=1pt] {$c'$}
 (0,-.5) circle (2pt) node [left=1pt] {$c''$}
 (1,0) circle (2pt) node [below=1pt] {$b_0$}
 (2,0) circle (2pt) node [below=1pt] {$b_1$}
 (4,0) circle (2pt) node [below=1pt] {$b_{n-5}$}
 (5,0) circle (2pt) node [below=1pt] {$b_{n-4}$}
 (6,0) circle (2pt) node [below=1pt] {$c$};
\end{tikzpicture}
}}
\]
In particular, $\hh^\Delta_a=1$ if $a\in\{c,c',c''\}$, and $\hh^\Delta_a=2$ otherwise.
Applying Lemmas~\ref{only1nonzero}(4) and~\ref{fundamental}(3), we get the following corollaries.

\begin{cor} \label{homboundby2}
$[\vv_{p,a},\vv_{q,b}]\leq 2$ and the inequality is strict if at least one of the vertices $a$ and $b$
belongs to $\{c,c',c''\}$.
\qed
\end{cor}

\begin{cor} \label{deltaboundby2}
Let $\sigma\colon A\to B\to C\to A[1]$ be a triangle in $\CD^b(Q)$ with $A$ or $C$ indecomposable.
Then $\delta_\sigma(\mm_{p,a})\leq 2$ and the inequality is strict if $a$ belongs to $\{c,c',c''\}$.
\qed
\end{cor}

The main aim of this subsection is to prove the following fact.

\begin{prop} \label{mainproptriang}
Let $\sigma\colon A\to B\to C\xrightarrow{\gamma}A[1]$ be a triangle in $\CD^b(Q)$ with $A$ and $C$ indecomposable.
Assume that $M,N\in\CD^b(Q)$ satisfy $\dimv M=\dimv N$, $\delta_{M,N}\geq 0$, $\supp(\delta_\sigma)\subseteq\supp(\delta_{M,N})$,
but the inequality $\delta_{M,N}\geq\delta_\sigma$ does not hold.

Then there is an indecomposable direct summand $C'$ of $N$ together with a morphism $h\colon C'\to C$
such that the pullback $\sigma'\colon A\to B'\to C'\xrightarrow{\gamma h}A[1]$ of $\sigma$ along $h$ does not split,
$\delta_{\sigma'}\leq\delta_{M,N}$ and $\supp(\delta_\sigma-\delta_{\sigma'})\subseteq\supp(\delta_{M,N}-\delta_{\sigma'})$.
\end{prop}

We introduce two integer-valued functions $\varphi$ and $\psi$ on $(\BZ\Delta)_0\,\sqcup\,(\BZ\Delta)_2$
as the compositions of the canonical bijection $(\BZ\Delta)_0\,\sqcup\,(\BZ\Delta)_2\to\BZ\times\Delta_0$
followed by the maps
\[
(p,a)\mapsto p+\dist(c,a)\quad\text{and}\quad(p,a)\mapsto p-\dist(c,a),
\]
respectively.

\begin{lem} \label{ifmesh2}
Let $\sigma\colon A\to B\to C\to A[1]$ be a triangle in $\CD^b(Q)$ such that $A$ and $C$ are indecomposable.
Let $\mm$ be a mesh in $\BZ\Delta$ such that $\delta_\sigma(\mm)=2$
\textup{(}in particular, $\mm=\mm_{p_0,b_r}$ for some pair $(p_0,b_r)$ in $(\BZ\times\Delta_0)^\mm$\textup{)}.
Then
\[
\delta_\sigma(\mm_{p,a})=\hh^\Delta_a
\]
for all meshes $\mm_{p,a}$ satisfying $\varphi(\mm_{p,a})\geq\varphi(\mm)$ and $\psi(\mm_{p,a})\leq\psi(\mm)$.
\end{lem}

We can illustrate the above statement about the function $\delta_\sigma$ for $r=2$ (hence $n\geq 6$)
by the following picture
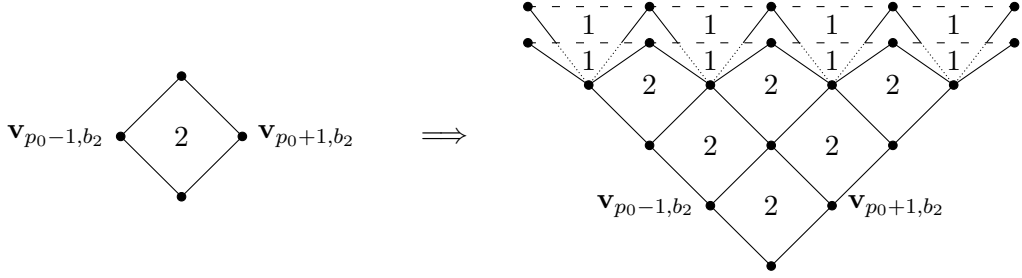
\[
\vcenter{\hbox{
\begin{tikzpicture} [>=stealth,scale=.8]
\filldraw
 (4,2) circle (2pt)
 (3,1) circle (2pt) node [left=2pt] {$\vv_{p_0-1,b_2}$}
 (5,1) circle (2pt) node [right=2pt] {$\vv_{p_0+1,b_2}$}
 (4,0) circle (2pt);
\draw
 node at (4,1) {$2$};
\draw (3,1) -- (4,2) -- (5,1);
\draw (3,1) -- (4,0) -- (5,1);
\end{tikzpicture}
}}
\quad\implies\quad
\vcenter{\hbox{
\begin{tikzpicture} [>=stealth,scale=.8]
\filldraw
 \foreach \x in {0,2,4,6,8}
  {
   (\x,4.3) circle (2pt)
   (\x,3.7) circle (2pt)
  }
 (1,3) circle (2pt)
 (3,3) circle (2pt)
 (5,3) circle (2pt)
 (7,3) circle (2pt)
 (2,2) circle (2pt)
 (4,2) circle (2pt)
 (6,2) circle (2pt)
 (3,1) circle (2pt) node [left=2pt] {$\vv_{p_0-1,b_2}$}
 (5,1) circle (2pt) node [right=2pt] {$\vv_{p_0+1,b_2}$}
 (4,0) circle (2pt);
\draw
 \foreach \x in {1,3,5,7}
  {
   node at (\x,4) {$1$}
   node at (\x,3.4) {$1$}
  }
 \foreach \x in {2,4,6} {node at (\x,3) {$2$}}
 \foreach \x in {3,5} {node at (\x,2) {$2$}}
 node at (4,1) {$2$};
\draw [loosely dashed] (0,4.3) -- (8,4.3);
\draw [loosely dashed] (0,3.7) -- (8,3.7);
\foreach \x in {0,2,4,6}
 {
  \draw (\x,4.3) -- (\x+6/13,3.7);
  \draw [densely dotted] (\x+6/13,3.7) -- (\x+1,3) -- (\x+20/13,3.7);
  \draw (\x+20/13,3.7) -- (\x+2,4.3);
  \draw (\x,3.7) -- (\x+1,3) -- (\x+2,3.7);
 }
\foreach \x in {1,3,5} \draw (\x,3) -- (\x+1,2) -- (\x+2,3);
\foreach \x in {2,4} \draw (\x,2) -- (\x+1,1) -- (\x+2,2);
\draw (3,1) -- (4,0) -- (5,1);
\end{tikzpicture}
}}
\]

\begin{proof}
The claim follows by induction on $r\geq 0$ from the following two properties of $\delta_\sigma$:
\begin{enumerate}
\item[(i)] If $r>0$ then $\delta_\sigma(\mm_{p_0-1,b_{r-1}})=\delta_\sigma(\mm_{p_0+1,b_{r-1}})=2$.
\item[(ii)] If $r=0$ then $\delta_\sigma(\mm_{p_0-1,c'})=\delta_\sigma(\mm_{p_0-1,c''})=\delta_\sigma(\mm_{p_0+1,c'})=\delta_\sigma(\mm_{p_0+1,c''})=1$.
\end{enumerate}
Indeed, the base step follows from (ii) and the induction step from (i).

Combining the assumption $\delta_\sigma(\mm_{p_0,b_r})=2$ with Lemma~\ref{fundamental}(3) and Corollary~\ref{homboundby2} gives the equalities
$[A,\vv_{p_0-1,b_r}]=2=[\vv_{p_0+1,b_r},C]$.
Since $A$ and $C$ are indecomposable, $[\vv_{p_0+1,b_r},A]=0=[C,\vv_{p_0-1,b_r}]$,
by Lemma~\ref{homifpath}(3).

Let $\CL$ be the set of vertices lying on the following sectional path in $\BZ\Delta$:
\[
\vv_{p_0-n+r+2,c}\to\vv_{p_0-n+r+3,b_{n-4}}\to\quad\cdots\quad\to\vv_{p_0-1,b_r}.
\]
Let $p'=p_0-n+r+1$.
Applying Corollary~\ref{multMNsection} for $\vv\in\CL$ we get
\[
\sum_{\vv\in\CL}\bigl(\mult_\vv(A\oplus C)-\mult_\vv(B)\bigr)=\delta_\sigma(\mm_{p',c})+\delta_\sigma(\mm_{p_0,b_r})
-\begin{cases}
\delta_\sigma(\mm_{p_0-1,b_{r-1}}),&r>0,\\
\delta_\sigma(\mm_{p_0-1,c'})+\delta_\sigma(\mm_{p_0-1,c''}),&r=0.
\end{cases}
\]
For any $\vv\in\CL$, $[\vv,\vv_{p_0-1,b_r}]=1$, by Lemma~\ref{hom1section}, thus $\vv$ is isomorphic neither to $A$ nor to $C$,
and consequently $\mult_\vv(A\oplus C)=0$.
Remembering that $\delta_\sigma(\mm_{p_0,b_r})=2$, we get
\[
2\leq
\begin{cases}
\delta_\sigma(\mm_{p_0-1,b_{r-1}}),&r>0,\\
\delta_\sigma(\mm_{p_0-1,c'})+\delta_\sigma(\mm_{p_0-1,c''}),&r=0.
\end{cases}
\]
On the other hand, by Corollary~\ref{deltaboundby2},
$\delta_\sigma(\mm_{p_0-1,b_{r-1}})\leq 2$, $\delta_\sigma(\mm_{p_0-1,c'})\leq 1$ and $\delta_\sigma(\mm_{p_0-1,c''})\leq 1$.
Consequently, we get three equalities of the six equalities appearing in (i) and (ii).

Dual considerations for the following sectional path in $\BZ\Delta$:
\[
\vv_{p_0+1,b_r}\to\quad\cdots\quad\to\vv_{p_0+n-r-3,b_{n-4}}\to\vv_{p_0+n-r-2,c}
\]
lead to the remaining three equalities in (i) and (ii).
\end{proof}

\begin{proof}[Proof of Proposition~\ref{mainproptriang}.]
Let $\mm_{p_0,a_0}$ be a mesh satisfying $\delta_{M,N}(\mm_{p_0,a_0})<\delta_\sigma(\mm_{p_0,a_0})$.
By Corollary~\ref{finitesupport}, we can choose such a mesh having the minimal value $\psi(\mm_{p_0,a_0})$.
We conclude from the assumption $\supp(\delta_\sigma)\subseteq\supp(\delta_{M,N})$ and Corollary~\ref{deltaboundby2} that
$\delta_{M,N}(\mm_{p_0,a_0})=1$, $\delta_\sigma(\mm_{p_0,a_0})=2$ and $a_0=b_r$ for some $r$ with $0\leq r\leq n-4$.
Let $\CR$ be the subset of $(\BZ\Delta)_0$ consisting of vertices $\vv_{p,a}$ such that
$\varphi(\vv_{p,a})>\varphi(\mm_{p_0,b_r})$ and $\psi(\vv_{p,a})<\psi(\mm_{p_0,b_r})$.
We illustrate the set $\CR$ for $r=3$ (hence $n\geq 7$) using $14$ big dots:
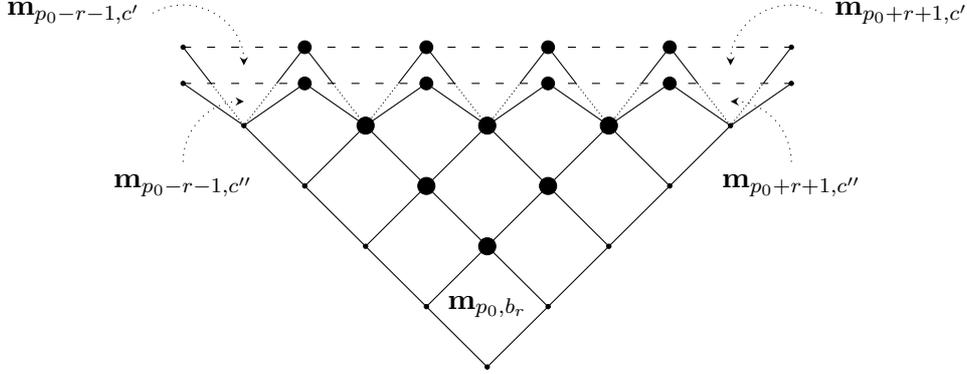
\[
\vcenter{\hbox{
\begin{tikzpicture} [>=stealth,scale=.8]
\filldraw
 (0,5.3) circle (1pt)
 (0,4.7) circle (1pt)
 \foreach \x in {2,4,6,8}
  {
   (\x,5.3) circle (3pt)
   (\x,4.7) circle (3pt)
  }
 (10,5.3) circle (1pt)
 (10,4.7) circle (1pt)
 (1,4) circle (1pt)
 (3,4) circle (4pt)
 (5,4) circle (4pt)
 (7,4) circle (4pt)
 (9,4) circle (1pt)
 (2,3) circle (1pt)
 (4,3) circle (4pt)
 (6,3) circle (4pt)
 (8,3) circle (1pt)
 (3,2) circle (1pt)
 (5,2) circle (4pt)
 (7,2) circle (1pt)
 (4,1) circle (1pt)
 (6,1) circle (1pt)
 (5,0) circle (1pt);
\draw node at (5,1) {$\mm_{p_0,b_r}$};
\draw (1,5) [<-,dotted] arc (0:120:1) node [left=1pt] {$\mm_{p_0-r-1,c'}$};
\draw (1,4.4) [<-,dotted] arc (90:180:1) node [below=1pt] {$\mm_{p_0-r-1,c''}$};
\draw (9,5) [<-,dotted] arc (180:60:1) node [right=1pt] {$\mm_{p_0+r+1,c'}$};
\draw (9,4.4) [<-,dotted] arc (90:0:1) node [below=1pt] {$\mm_{p_0+r+1,c''}$};
\draw [loosely dashed] (0,5.3) -- (10,5.3);
\draw [loosely dashed] (0,4.7) -- (10,4.7);
\foreach \x in {0,2,4,6,8}
 {
  \draw (\x,5.3) -- (\x+6/13,4.7);
  \draw [densely dotted] (\x+6/13,4.7) -- (\x+1,4) -- (\x+20/13,4.7);
  \draw (\x+20/13,4.7) -- (\x+2,5.3);
  \draw (\x,4.7) -- (\x+1,4) -- (\x+2,4.7);
 }
\foreach \x in {1,3,5,7} \draw (\x,4) -- (\x+1,3) -- (\x+2,4);
\foreach \x in {2,4,6} \draw (\x,3) -- (\x+1,2) -- (\x+2,3);
\foreach \x in {3,5} \draw (\x,2) -- (\x+1,1) -- (\x+2,2);
\draw (4,1) -- (5,0) -- (6,1);
\end{tikzpicture}
}}
\]
The key observation is that by Corollary~\ref{multMN} we get the formula
\begin{multline*}
\sum_{\vv_{p,a}\in\CR}\,\hh^\Delta_a\cdot\bigl(\mult_{\vv_{p,a}}(N)-\mult_{\vv_{p,a}}(M)\bigr)
=-2\cdot\delta_{M,N}(\mm_{p_0,b_r})\\
+\delta_{M,N}(\mm_{p_0-r-1,c'})+\delta_{M,N}(\mm_{p_0-r-1,c''})+\delta_{M,N}(\mm_{p_0+r+1,c'})+\delta_{M,N}(\mm_{p_0+r+1,c''}).
\end{multline*}
Combining Lemma~\ref{ifmesh2} with the assumption $\supp(\delta_\sigma)\subseteq\supp(\delta_{M,N})$
and using $\delta_{M,N}(\mm_{p_0,b_r})=1$, we get that the right-hand side is at least $2$.
Hence $\mult_{C'}(N)>0$ for some $C'=\vv_{p',a'}$ in $\CR$.
Again by Lemma~\ref{ifmesh2}, $\delta_\sigma(\mm_{p'-1,a'})>0$, which from the definition of $\delta_\sigma$ means that
$\gamma\circ h\neq 0$ for some morphism $h\colon C'\to C$.
Let $\sigma'$ be the pullback of $\sigma$ along $h$.
We need to prove that $\delta_{\sigma'}\leq\delta_{M,N}$ and $\supp(\delta_\sigma-\delta_{\sigma'})\subseteq\supp(\delta_{M,N}-\delta_{\sigma'})$.

Since $\delta_{\sigma'}\leq\delta_\sigma$ and $\supp(\delta_\sigma)\subseteq\supp(\delta_{M,N})$, it suffices to show that
$\delta_{M,N}(\mm_{p,a})\geq\delta_\sigma(\mm_{p,a})$ whenever $\delta_{\sigma'}(\mm_{p,a})>0$.
Thus we assume that $\delta_{\sigma'}(\mm_{p,a})>0$.
By Lemma~\ref{fundamental}(3), $[\vv_{p+1,a},C']>0$, and from Lemma~\ref{homifpath}(3) we conclude that $\psi(\vv_{p+1,a})\leq\psi(C')$.
Using the fact that $C'$ belongs to $\CR$ and how the latter was defined, we get the following sequence of inequalities
\[
\psi(\mm_{p,a})<\psi(\vv_{p+1,a})\leq\psi(C')<\psi(\mm_{p_0,b_r}).
\]
It follows from our choice of the mesh $\mm_{p_0,b_r}$ that $\delta_{M,N}(\mm_{p,a})\geq\delta_\sigma(\mm_{p,a})$,
which finishes the proof.
\end{proof}

\subsection{Passage from $\CD^b(Q)$ to $\rep(Q)$}

The main aim of this subsection is to prove a result analogous to Proposition~\ref{mainproptriang},
concerning the category $\rep(Q)$, where $Q$ is a Dynkin quiver of type $\BD$.
Throughout this subsection $Q$ is a Dynkin quiver with its underlying graph $\Delta$.

As observed in Subsection~\ref{subsect mesh}, the $\Bbbk$-linear structure of the category $\CD^b(Q)$
is fully described by the quiver $\BZ\Delta$.
Similarly, the category $\rep(Q)$ is fully described by its Auslander-Reiten quiver $\Gamma_Q$.
Moreover, the identification of $\rep(Q)$ as a full subcategory of $\CD^b(Q)$ correspond to
the identification of $\Gamma_Q$ as a full convex subquiver of $\BZ\Delta$, which we are going to explain.
We note that introducing the Auslander-Reiten quiver $\Gamma_Q$ as a subquiver of $\BZ\Delta$ was done
already in~\cite{Gab}*{6.5}.

By a slice in $\BZ\Delta$ we mean a full convex subquiver containing exactly
one vertex $\vv_{r_a,a}$ for each $a\in\Delta_0$.
Thus $|r_a-r_b|=1$ for any adjacent vertices $a$ and $b$.
Recall that $P_a=\vv_{p_a,a}$ for any vertex $a\in Q_0=\Delta_0$.
The vertices $P_a$, $a\in Q_0$, together with the arrows connecting them form a slice $\CS$ isomorphic to $Q^{\op}$,
where $Q^{\op}$ is the opposite quiver of $Q$ having the same set of vertices, but with the arrows reversed.
Consequently, the vertices $I_a=\nu(P_a)=\vv_{p_a+n_\Delta-2,\phi_\Delta(a)}$, $a\in Q_0$, lie on a slice $\nu\CS$,
which is also isomorphic to $Q^{\op}$.
Then the Auslander-Reiten quiver $\Gamma_Q$ is the smallest full convex subquiver of $\BZ\Delta$
containing $\CS$ and $\nu\CS$.
We denote by $(\Gamma_Q)_2$ the set of all meshes in $\BZ\Delta$ which are contained in $\Gamma_Q$.

The shifts $(\Gamma_Q)[i]$, $i\in\BZ$, are pairwise disjoint subquivers of $\BZ\Delta$, hence we have the following inclusion
\[
\bigsqcup_{i\in\BZ}\,(\Gamma_Q)[i]\,\subseteq\,\BZ\Delta.
\]
In fact, this inclusion is the equality on the sets of vertices, and only the arrows connecting $\nu\CS[i]$ with $\CS[i+1]$, $i\in\BZ$, are missing
(see for instance~\cite{Hap}*{5.5}).
For example, if $Q$ is the quiver
\[
\begin{tikzpicture} [>=stealth,shorten <= 6pt,shorten >=6pt]
\filldraw
 (0,0) circle (2pt) node [above=1pt] {$a$}
 (1,0) circle (2pt) node [above=1pt] {$b$}
 (2,0) circle (2pt) node [above=1pt] {$c$}
 (3,0) circle (2pt) node [above=1pt] {$d$}
;
\draw [->] (0,0) -- (1,0);
\draw [->] (2,0) -- (1,0);
\draw [->] (3,0) -- (2,0);
\end{tikzpicture}
\]
then the above embedding of quivers looks as follows.
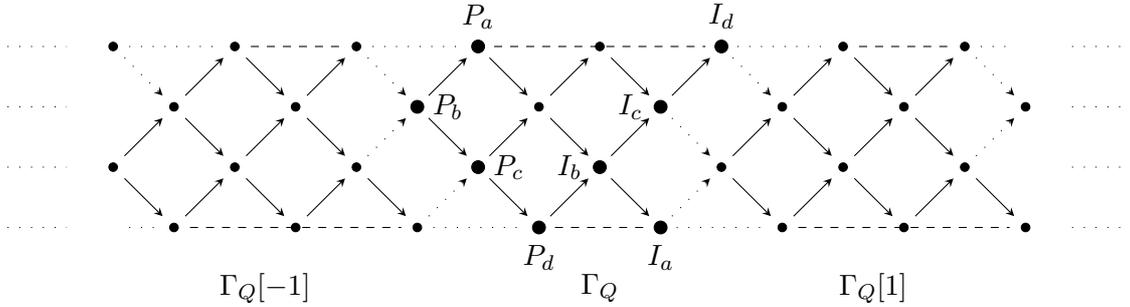
\[
\begin{tikzpicture} [>=stealth,shorten <= 6pt,shorten >=6pt,scale=.8]
\filldraw
 \foreach \x in {0,2,...,14}
  {
   (\x,3) circle (2pt)
   (\x+1,2) circle (2pt)
   (\x,1) circle (2pt)
   (\x+1,0) circle (2pt)
  }
 (6,3) circle (3pt)
 (10,3) circle (3pt)
 (5,2) circle (3pt)
 (9,2) circle (3pt)
 (6,1) circle (3pt)
 (8,1) circle (3pt)
 (7,0) circle (3pt)
 (9,0) circle (3pt)
 ;
\foreach \x in {1,3,...,13}
 {
  \draw [->] (\x,2) -- (\x+1,3);
 }
\foreach \x in {2,6,8,12}
 {
  \draw [->] (\x,3) -- (\x+1,2);
 }
\foreach \x in {0,2,6,8,10,12}
 {
  \draw [->] (\x,1) -- (\x+1,2);
 }
\foreach \x in {1,3,5,7,11,13}
 {
  \draw [->] (\x,2) -- (\x+1,1);
 }
\foreach \x in {0,2,...,14}
 {
  \draw [->] (\x,1) -- (\x+1,0);
 }
\foreach \x in {1,3,7,11,13}
 {
  \draw [->] (\x,0) -- (\x+1,1);
 }
\draw [->,loosely dotted] (0,3) -- (1,2);
\draw [->,loosely dotted] (4,3) -- (5,2);
\draw [->,loosely dotted] (10,3) -- (11,2);
\draw [->,loosely dotted] (14,3) -- (15,2);
\draw [->,loosely dotted] (4,1) -- (5,2);
\draw [->,loosely dotted] (9,2) -- (10,1);
\draw [->,loosely dotted] (14,1) -- (15,2);
\draw [->,loosely dotted] (5,0) -- (6,1);
\draw [->,loosely dotted] (9,0) -- (10,1);
\node at (2.5,-1) {$\Gamma_Q[-1]$};
\node at (8,-1) {$\Gamma_Q$};
\node at (12.5,-1) {$\Gamma_Q[1]$};
\node at (6,3.5) {$P_a$};
\node at (10,3.5) {$I_d$};
\node at (5.5,2) {$P_b$};
\node at (8.5,2) {$I_c$};
\node at (6.5,1) {$P_c$};
\node at (7.5,1) {$I_b$};
\node at (7,-.5) {$P_d$};
\node at (9,-.5) {$I_a$};
\draw [dashed] (2,3) -- (4,3);
\draw [dashed] (6,3) -- (10,3);
\draw [dashed] (12,3) -- (14,3);
\draw [dashed] (1,0) -- (5,0);
\draw [dashed] (7,0) -- (9,0);
\draw [dashed] (11,0) -- (15,0);
\draw [loosely dotted] (0,3) -- (2,3);
\draw [loosely dotted] (4,3) -- (6,3);
\draw [loosely dotted] (10,3) -- (12,3);
\draw [loosely dotted] (14,3) -- (15,3);
\draw [loosely dotted] (0,0) -- (1,0);
\draw [loosely dotted] (5,0) -- (7,0);
\draw [loosely dotted] (9,0) -- (11,0);
\foreach \y in {0,1,2,3}
 {
  \draw [loosely dotted] (-2,\y) -- (-.5,\y);
  \draw [loosely dotted] (15.5,\y) -- (17,\y);
 }
\end{tikzpicture}
\]

Let $X$ and $Y$ be representations in $\rep(Q)$.
Then $[X,Y]^i=0$ for $i<0$, and hence $\langle X,Y\rangle=[X,Y]$.
Therefore the definition of the integer-valued function $\delta_{M,N}$ simplifies
if we restrict to the subcategory $\rep(Q)$ as the following result explains.

\begin{lem}
Let $Q$ be a Dynkin quiver.
Assume that $M$ and $N$ belong to $\rep(Q)$ and $\dimv M=\dimv N$.
Let $\mm_{p,a}$ be a mesh in $(\BZ\Delta)_2$.
Then
\[
\delta_{M,N}(\mm_{p,a})=[\vv_{p+1,a},N]-[\vv_{p+1,a},M]
=[N,\vv_{p-1,a}]-[M,\vv_{p-1,a}]
\]
if $\mm_{p,a}$ belongs to $(\Gamma_Q)_2$, and $\delta_{M,N}(\mm_{p,a})=0$, otherwise.
\qed
\end{lem}

Recall that for exact sequence $\sigma\colon 0\to A\xrightarrow{\alpha}B\xrightarrow{\beta}C\to 0$ in $\rep(Q)$
we have the corresponding triangle
$\hat{\sigma}\colon A\xrightarrow{\alpha}B\xrightarrow{\beta}C\xrightarrow{\gamma}A[1]$ in $\CD^b(Q)$.
We set $\delta_\sigma=\delta_{\hat{\sigma}}$.
In this case Definition~\ref{defdeltasigma} and Corollary~\ref{deltasigmaalter} simplify
as shown in the following lemma.
We also note that the function $\delta_\sigma$ is closely related to the defect functors considered in~\cite{ARS}*{IV.4}.

\begin{lem} \label{delsigmaforreps}
Let $\sigma\colon 0\to A\xrightarrow{\alpha}B\xrightarrow{\beta}C\to 0$ be a short exact sequence in $\rep(Q)$
and $\mm_{p,a}$ be any mesh in $\BZ\Delta$.
Then $\delta_{\sigma}(\mm_{p,a})=0$ if $\mm_{p,a}$ does not belong to $(\Gamma_Q)_2$.
Otherwise,
\[
\delta_{\sigma}(\mm_{p,a})=[\vv_{p+1,a},A\oplus C]-[\vv_{p+1,a},B]
=[A\oplus C,\vv_{p-1,a}]-[B,\vv_{p-1,a}]
\]
equals the ranks of the last $\Bbbk$-linear morphisms in the induced exact sequences:
\[
0\to\Hom_Q(\vv_{p+1,a},A)\to\Hom_Q(\vv_{p+1,a},B)\to\Hom_Q(\vv_{p+1,a},C)\to\Ext^1_Q(\vv_{p+1,a},A)
\]
and
\[
0\to\Hom_Q(C,\vv_{p-1,a})\to\Hom_Q(B,\vv_{p-1,a})\to\Hom_Q(A,\vv_{p-1,a})\to\Ext^1_Q(C,\vv_{p-1,a}). \eqno \qed
\]
\end{lem}

We are ready to restrict Proposition~\ref{mainproptriang} to the category $\rep(Q)$.

\begin{cor} \label{maincortriang}
Let $Q$ be a Dynkin quiver of type $\BD$ and $\sigma\colon 0\to A\to B\to C\to 0$
be a short exact sequence in $\rep(Q)$ with $A$ and $C$ indecomposable.
Assume that $M,N\in\rep(Q)$ satisfy $\dimv M=\dimv N$, $\delta_{M,N}\geq 0$, $\supp(\delta_\sigma)\subseteq\supp(\delta_{M,N})$,
but the inequality $\delta_{M,N}\geq\delta_\sigma$ does not hold.

Then there is an indecomposable direct summand $C'$ of $N$ together with a homomorphism $h\colon C'\to C$
such that the pullback $\sigma'\colon A\to B'\to C'\to 0$ of $\sigma$ along $h$ does not split,
$\delta_{\sigma'}\leq\delta_{M,N}$ and $\supp(\delta_\sigma-\delta_{\sigma'})\subseteq\supp(\delta_{M,N}-\delta_{\sigma'})$.
\qed
\end{cor}

\section{Schemes of representations of quivers}
\label{sectionproof}

The set $\rep_Q^\dd(R)$, where $R$ is a commutative $\Bbbk$-algebra, has a natural structure of a vector space over $\Bbbk$,
and using the addition simplifies working with elements of this set.
In particular, any element of $\rep_Q^\dd(\Bbbk[\varepsilon])$ can be uniquely presented in the form
$N+\varepsilon\cdot Z$, where $N$ and $Z$ belong to $\rep_Q^\dd(\Bbbk)$.
At the same time $Z$ is viewed as a tangent vector in the Zariski tangent space $T_N\rep_Q^\dd$.
Hence $T_N \rep_Q^\dd = \rep_Q^\dd (\Bbbk)$ for any $N \in \rep_Q^\dd (\Bbbk)$.
We start this section with a more adequate representation-theoretic interpretation of elements of $T_N \rep_Q^\dd$.

\subsection{Tangent vectors and short exact sequences}

Throughout this subsection $Q$ is a finite quiver.
If $U=(U_a;\,a\in Q_0)$ and $V=(V_a;\,a\in Q_0)$ are collections of finite dimensional vector spaces indexed by the vertices of $Q$,
then we define $\Bbbk$-schemes $\CV^U_V$ and $\CA^U_V$ by
\[
\CV^U_V(R)=\prod_{a\in Q_0}\Hom_R(R \otimes U_a, R \otimes V_a) \qquad \text{and} \qquad
\CA^U_V(R)=\prod_{\alpha \in Q_1} \Hom_R (R \otimes U_{s \alpha}, R \otimes V_{t \alpha}),
\]
where all tensor products in this section are taken over $\Bbbk$.
In particular, if $\Bbbk^{\dd}=(k^{\dd_a})$
and we identify the $R$-homomorphisms of the form $R \otimes k^d \to R \otimes k^e$ with the corresponding matrices in $\BM_{e \times d} (R)$, then
\begin{equation} \label{eq_repA}
\CA^{\Bbbk^\dd}_{\Bbbk^\dd}=\rep_Q^\dd.
\end{equation}
Similarly $\GL_\dd$ is an open subscheme of $\CV_{\Bbbk^\dd}^{\Bbbk^\dd}$.

If $h\in\CV^U_V(R)$ and $Z\in\CA^V_W(R)$, then we define $Z \circ h \in \CA^U_W(R)$ by
\[
(Z \circ h)_\alpha = Z_\alpha \circ h_{s \alpha}.
\]
for any  $\alpha \in Q_1$.
We define $\circ \colon \CV^V_W (R) \times \CA^U_V (R) \to \CA^U_W (R)$ dually.
Finally we have $\circ \colon \CV^V_W (R) \times \CV^U_V (R) \to \CV^U_W (R)$ defined in an obvious way.
We note that the action $\star$ of $\GL_\dd\subseteq\CV_{\Bbbk^\dd}^{\Bbbk^\dd}$ on $\rep_Q^\dd=\CA^{\Bbbk^\dd}_{\Bbbk^\dd}$
can be written as $g\star M=g\circ M\circ g^{-1}$.

If $U = \bigoplus_{s \in \CS} U^s$ (i.e.\ $U_a = \bigoplus_{s \in \CS} U_a^s$ for each $a \in Q_0$)
and $V = \bigoplus_{s' \in \CS'} V^{s'}$, then
\begin{equation} \label{eq_decomp}
\CV^U_V = \prod_{s \in \CS, \; s' \in \CS'} \CV^{U^s}_{V^{s'}}
\qquad \text{and} \qquad
\CA^U_V = \prod_{s \in \CS, \; s' \in \CS'} \CA^{U^s}_{V^{s'}}.
\end{equation}
For $h \in \CV^U_V (R)$, $q \in \CS$ and $p \in \CS'$,
we denote by $h^{p, q}$ the image of $h$ under the projection $\CV^U_V (R) \to \CV^{U^q}_{V^p} (R)$ induced by the former of the above decompositions.
Conversely, given $h^{p, q} \in \CV^{U^q}_{V^p} (R)$, then $\widehat{h^{p, q}}$ denotes the image of $h^{p, q}$ under the section
$\CV^{U^q}_{V^p} (R) \to \CV^U_V (R)$.
We define $Z^{p, q} \in \CA^{U^q}_{V^p} (R)$, for $Z \in \CA^U_V (R)$, and $\widehat{Z^{p, q}} \in \CA^U_V (R)$, for $Z^{p, q} \in \CA^{U^q}_{V^p} (R)$, analogously.
In particular,
\begin{equation} \label{eq_sum}
h=\sum_{s \in \CS, \; s' \in \CS'}\widehat{h^{s',s}} \qquad \text{and} \qquad
Z=\sum_{s \in \CS, \; s' \in \CS'}\widehat{Z^{s',s}}.
\end{equation}

If $U=(U_a,U_\alpha)$ and $V=(V_a,V_\alpha)$ are representations in $\rep (Q)$, then we put
\[
\CV^U_V = \CV^{(U_a)}_{(V_a)} \qquad \text{and} \qquad \CA^U_V = \CA^{(U_a)}_{(V_a)}.
\]
Moreover, we define $\BZ_Q^1 (U, V) = \CA^U_V (\Bbbk)$.
Since $\Hom_Q (U, V) \subseteq \CV^U_V (\Bbbk)$,
we obtain using the above defined compositions that $\BZ_Q^1 (?, ?)$ is a functor from $\rep(Q)^{\op}\times\rep(Q)$ to $\mod \Bbbk$.
Observe that if $N \in \rep_Q^\dd(\Bbbk)$, then
\begin{equation} \label{eq_tanA}
\BZ_Q^1 (N, N) = \CA^{\Bbbk^\dd}_{\Bbbk^\dd} (\Bbbk) = \rep_Q^\dd (\Bbbk) = T_N \rep_Q^\dd,
\end{equation}
and we will view the elements of $\BZ_Q^1 (N, N)$ as tangent vectors.

For representations $U$ and $V$ in $\rep(Q)$ and $Z \in \BZ_Q^1 (V, U)$,
let $W=W(U,Z,V)$ denote the representation such that
$W_a=U_a\oplus V_a$ and
\[
W_\alpha\colon \, U_{s\alpha}\oplus V_{s\alpha}\xrightarrow{\bsmatrix{U_\alpha&Z_\alpha\\ 0&V_\alpha}}U_{t\alpha}\oplus V_{t\alpha}.
\]
Let $\varphi\colon U\to W$ and $\psi \colon W\to V$ denote the homomorphisms such that $\varphi_a\colon U_a\to U_a\oplus V_a$
is the canonical section and $\psi_a\colon U_a\oplus V_a\to V_a$ is the canonical projection, for any $a\in Q_0$.
Then
\[
\sigma(U,Z,V)\colon\,0\to U\xrightarrow{\varphi}W(U,Z,V)\xrightarrow{\psi}V\to 0
\]
is a short exact sequence in $\rep(Q)$.
If it causes no confusion, we will present symbolically the sequence $\sigma(U,Z,V)$ in the form
\[
0\to U\xrightarrow{\bsmatrix{1\\ 0}}\bsmatrix{U&Z\\ 0&V}\xrightarrow{\bsmatrix{0&1}}V\to 0.
\]
Moreover, if $h\in\Hom_Q(V',V)$ then $\sigma(U,Z \circ h,V')$ is the pullback of $\sigma(U,Z,V)$ along $h$
leading to the following commutative diagram
\[
\vcenter{\hbox{
\xymatrix{
0\ar[r]&U\ar[r]^-{\bsmatrix{1\\ 0}}\ar@{=}[d]&{\bsmatrix{U&Z \circ h\\ 0&V'}}\ar[r]^-{\bsmatrix{0&1}}\ar[d]^-{\bsmatrix{1&0\\ 0&h}}&V'\ar[r]\ar[d]^-h&0\\
0\ar[r]&U\ar[r]^-{\bsmatrix{1\\ 0}}&{\bsmatrix{U&Z\\ 0&V}}\ar[r]^-{\bsmatrix{0&1}}&V\ar[r]&0
}}}
\]
Dually $\sigma(U',h' \circ Z,V)$ is the pushout of $\sigma(U,Z,V)$ along $h'$, for any homomorphism $h'\colon U\to U'$.

Let $U$ and $V$ be representations in $\rep(Q)$. We denote by $\BB^1_Q(V,U)$ the image of the map
\[
\eta_{V,U}\colon\CV^V_U(\Bbbk)\to\CA^V_U(\Bbbk)=\BZ_Q^1(V,U),
\qquad
\eta_{V,U}(h)=h \circ (V_\alpha) - (U_\alpha) \circ h.
\]
Then $\BB^1_Q(?,?)$ is a $\Bbbk$-linear subfunctor of $\BZ^1_Q(?,?)$.
Moreover, $Z\in\BZ^1_Q(V,U)$ belongs to $\BB^1_Q(V,U)$ if and only if the short exact sequence $\sigma(U,Z,V)$ splits.
The quotient of $\BZ^1_Q(V,U)$ by $\BB^1_Q(V,U)$ can be identified with the extension group $\Ext^1_Q(V,U)$ of $V$ by $U$,
and we have the following exact sequence
\[
0\to\Hom_Q(V,U)\to\CV^V_U(\Bbbk)\xrightarrow{\eta_{V,U}}\CA^V_U(\Bbbk)\to\Ext^1_Q(V,U)\to 0.
\]

Applying the above for $U = N = V$, where $N\in\rep_Q^\dd(\Bbbk)$, we get the exact sequence
\[
0\to\End_Q(N)\to\CV^{\Bbbk^\dd}_{\Bbbk^\dd}(\Bbbk)\xrightarrow{\eta_{N,N}}\CA^{\Bbbk^\dd}_{\Bbbk^\dd}(\Bbbk)\to\Ext^1_Q(N,N)\to 0.
\]
The space $\CV^{\Bbbk^\dd}_{\Bbbk^\dd}(\Bbbk)$ can be identified with the tangent space $T_1\GL_\dd$,
the space $\CA^{\Bbbk^\dd}_{\Bbbk^\dd}(\Bbbk)=\BZ^1_Q(N,N)$ with the tangent space $T_N\rep_Q^\dd$,
and $\eta_{N,N}$ with the tangent map induced by the orbit map
\[
\GL_\dd\to\rep_Q^\dd,\quad g\mapsto g\star N.
\]
Under this identification, $T_N\CO_N=\BB^1_Q(N,N)$ and the normal space $T_N\rep_Q^\dd/T_N\CO_N$
at $N$ to $\CO_N$ in $\rep_Q^\dd$ coincides with $\Ext^1_Q(N,N)$, which is a famous Voigt result~\cite{Voi}.

\subsection{Proof of the main result}

Throughout this subsection $Q$ is a Dynkin quiver, $M$ a representation in $\rep(Q)$ and $\dd=\dimv M$.
We will work with closed subschemes of $\rep_Q^\dd$ and hence we start with a few general remarks on subschemes.
If $\CX$ is a subscheme of a $\Bbbk$-scheme $\CY$, then the corresponding map $\CX(R)\to\CY(R)$ is injective,
and we will identify $\CX(R)$ with its image in $\CY(R)$, for any commutative $\Bbbk$-algebra $R$.
The following fact can be concluded from~\cite{DG}*{I.2.6.1}:

\begin{lem} \label{propertyclosedsub}
Let $\CX$ be a closed subscheme of a $\Bbbk$-scheme $\CY$ and $\varphi\colon R\to S$ an injective homomorphism
of commutative $\Bbbk$-algebras.
Then $\CX(R)=\CY(\varphi)^{-1}(\CX(\CS))$. \qed
\end{lem}

We will use several times the above lemma for a $\Bbbk$-scheme $\CY$ which is affine (specifically for $\CY=\rep_Q^\dd$).
Then the closed embedding $\CX\subseteq\CY$ is isomorphic to
$\Spec(\psi)\colon\Spec(A/I)\to\Spec(A)$,
where $\psi\colon A\to A/I$ is the canonical surjective homomorphism,
for some commutative $\Bbbk$-algebra $A$ and its ideal $I$.
Hence the claim translates to an obvious fact about the existence of a homomorphism completing
a given commutative diagram in the category of $\Bbbk$-algebras
to another commutative diagram, as follows:
\[
\vcenter{\hbox{\xymatrix{
A\ar[r]\ar[d]_-\psi&R\ar[d]^-\varphi\\
A/I\ar[r]&S
}}}
\quad\implies\quad
\vcenter{\hbox{\xymatrix{
A\ar[r]\ar[d]_-\psi&R\ar[d]^-\varphi\\
A/I\ar@{-->}[ru]^-\exists\ar[r]&S
}}}
\]

If $\varphi\colon R\to S$ is an injective $\Bbbk$-algebra homomorphism, $\CX$ is a $\Bbbk$-scheme and $x\in\CX(R)$,
it will be convenient to denote the image of $x$ under the map $\CX(\varphi)\colon\CX(R)\to\CX(S)$ also by $x$.

Let $\pi\colon\Bbbk[\varepsilon]\to\Bbbk$ denote the canonical surjective homomorphism.
Since $\ov{\CO}_M$ is a closed subscheme of $\CC_M$ and the latter is a closed subscheme of $\rep_Q^\dd$,
we have the following commutative diagram with inclusions:
\[
\xymatrix{
\ov{\CO}_M(\Bbbk[\varepsilon])\ar@{}[r]|-{\subseteq}\ar[d]^-{\ov{\CO}_M(\pi)}
 &\CC_M(\Bbbk[\varepsilon])\ar@{}[r]|-{\subseteq}\ar[d]^{\CC_M(\pi)}
 &\rep_Q^\dd(\Bbbk[\varepsilon])\ar[d]^-{\rep_Q^\dd(\pi)}\\
\ov{\CO}_M(\Bbbk)\ar@{}[r]|-{=}&\CC_M(\Bbbk)\ar@{}[r]|-{\subseteq}&\rep_Q^\dd(\Bbbk)
}
\]
where the equality in the second row follows from (a reformulation of) Theorem~\ref{degeqhom}.
Observe that $\ov{\CO}_M(\Bbbk[\varepsilon])=\CC_M(\Bbbk[\varepsilon])$ if and only if
$T_N\ov{\CO}_M=T_N\CC_M$ for all points $N\in\ov{\CO}_M(\Bbbk)$.
It follows from~\cite{RZrank} that
\begin{align}
\label{descriptionCMone}
\ov{\CO}_M(\Bbbk)=\CC_M(\Bbbk)=\{N\in\rep_Q^\dd(\Bbbk);&\;\delta_{M,N}\geq 0\},\\
\label{descriptionCMtwo}
\CC_M(\Bbbk[\varepsilon])=\{N+\varepsilon\cdot Z\in\rep_Q^\dd(\Bbbk[\varepsilon]);&\;\text{$\delta_{M,N}\geq 0$ and
 $\supp(\delta_{\sigma(N,Z,N)})\subseteq\supp(\delta_{M,N})$}\}.
\end{align}
Let $N\in\CC_M(\Bbbk)=\ov{\CO}_M(\Bbbk)$.
Given $U$ and $V$ in $\rep(Q)$ we denote by $\BZ^1_{M,N}(V,U)$ the subset of $\BZ^1_Q(V,U)$ consisting
of the elements $Z$ such that $\supp(\delta_{\sigma(U,Z,V)})\subseteq\supp(\delta_{M,N})$.
In particular,
\[
T_N\CC_M=\BZ^1_{M,N}(N,N).
\]
It follows from~\cite{RZrank}*{Section 7} that $\BZ^1_{M,N}(?,?)$ is a $\Bbbk$-linear subfunctor of $\BZ^1_Q(?,?)$
containing $\BB^1_Q(?,?)$.

Let $N\in\rep_Q^\dd(\Bbbk)$.
Assume that we have a fixed decomposition $N=\bigoplus_{s\in\CS}N^s$ as a representation of $Q$.
In particular, the collection $(N_a)=\Bbbk^{\dd}$
decomposes as $\bigoplus_{s \in \CS} (N_a^s)$.
If $p, q \in \CS$,
then using notation introduced after~\eqref{eq_decomp} we have $H^{p,q} \in \CA^{N^q}_{N^p}(\Bbbk)$ for each $H \in \CA^{\Bbbk^\dd}_{\Bbbk^\dd}(\Bbbk)$.
Equality~\eqref{eq_tanA} implies that we can apply the above notation both when $H=L$ is a point of $\rep_Q^\dd(\Bbbk)$,
and when $H=Z$ is viewed as tangent vector.
Conversely, given $H^{p,q}\in\CA^{N^q}_{N^p}(\Bbbk)$ we have $\widehat{H^{p,q}}\in\CA^{\Bbbk^\dd}_{\Bbbk^\dd}(\Bbbk)$.
Observe that $N^{p,p}$ is the collection $(N_\alpha^p)$ and
$N^{p,q}=0$, for all $p \neq q$ in $\CS$,
hence $N=\sum_{s \in \CS} \widehat{N^{s,s}}$.
Analogously, applying notation introduced after~\eqref{eq_decomp} for the scheme $\CV^{\Bbbk^\dd}_{\Bbbk^\dd}(\Bbbk)$
we also have $g^{p, q} \in \CV^{N^q}_{N^p}(R)$ for $g\in\GL_\dd(R)$.
Finally, if $s\in\CS$, we set $1^s$ for the identity on $N^s$,
which is an element of $\CV_{N^s}^{N^s}(\Bbbk)$.
Consequently, $\sum_{s \in \CS} \widehat{1^s}$ is the identity on $(N_a)=\Bbbk^{\dd}$.

\begin{lem} \label{deltasigmapq}
Let $N$ be a point in $\rep_Q^\dd(\Bbbk)$ with a fixed decomposition $N=\bigoplus_{s\in\CS}N^s$.
Consider the representation $L=N+\widehat{Z^{p,q}}$ in $\rep_Q^\dd(\Bbbk)$
for $Z^{p,q} \in \CA^{N^q}_{N^p} (\Bbbk)=\BZ^1_Q(N^q,N^p)$ and two indices $p\neq q$ in $\CS$.
Then the following conditions hold:
\begin{enumerate}
\item[\textup{(1)}] $\delta_{L,N}=\delta_{\sigma(N^p,Z^{p,q},N^q)}$.
\item[\textup{(2)}] $L$ degenerates to $N$, i.e.\ $N\in\ov{\CO}_L(\Bbbk)$.
\item[\textup{(3)}] $L\simeq N$ if and only if $Z^{p,q}$ belongs to $\BB^1_Q(N^q,N^p)$.
\end{enumerate}
\end{lem}

\begin{proof}
Since $L=\sum_{s \in \CS} \widehat{N^{s,s}} + \widehat{Z^{p, q}}$, we have the following isomorphism of representations
\[
L\simeq\bigoplus_{s\neq p,q}N^s\oplus\bsmatrix{N^p&Z^{p,q}\\ 0&N^q},
\]
and (1) follows. In particular, $\delta_{L,N}\geq 0$, hence (2) holds, by \eqref{descriptionCMone}.

(3).
We have from (1) that $L$ is isomorphic to $N$ if and only if the sequence $\sigma(N^p,Z^{p,q},N^q)$
splits.
The latter means that $Z^{p,q}$ belongs to $\BB^1_Q(N^q,N^p)$.
\end{proof}

\begin{lem} \label{constrtangent}
Let $N=\bigoplus_{s\in\CS}N^s$ be a fixed decomposition of a point $N$ in $\ov{\CO}_M(\Bbbk)$.
Let $Z^{p,q}\in\CA^{N^q}_{N^p} (\Bbbk)$, for $p\neq q$ in $\CS$, be such that $\delta_{\sigma(N^p,Z^{p,q},N^q)}\leq\delta_{M,N}$.
Then
\[
N+\widehat{Z^{p,q}}\in\ov{\CO}_M(\Bbbk)
\qquad\text{and}\qquad
N+\varepsilon\cdot\widehat{Z^{p,q}}\in\ov{\CO}_M(\Bbbk[\varepsilon])
\qquad (\text{equivalently, }
\widehat{Z^{p,q}}\in T_N\ov{\CO}_M).
\]
\end{lem}

\begin{proof}
Let $L=N+\widehat{Z^{p,q}}$. By Lemma~\ref{deltasigmapq}(1), $\delta_{L,N}=\delta_{\sigma(N^p,Z^{p,q},N^q)}$ and
$\delta_{M,L}=\delta_{M,N}-\delta_{\sigma(N^p,Z^{p,q},N^q)}\geq 0$.
Hence $L\in\ov{\CO}_M(\Bbbk)$, by \eqref{descriptionCMone}.

We claim that the point $N+t\cdot\widehat{Z^{p,q}}$ of $\rep_Q^\dd(\Bbbk[t])$
belongs to $\ov{\CO}_M(\Bbbk[t])$.
Consider the element $g$ in $\GL_\dd(\Bbbk[t,t^{-1}])$ given by
$g=\sum_{s\neq p} \widehat{1^s}+t \cdot \widehat{1^p}$.
Obviously $g^{-1}=\sum_{s\neq p} \widehat{1^s}+t^{-1} \cdot \widehat{1^p}$.
Since $L=\sum_{s \in \CS} \widehat{N^{s,s}} + \widehat{Z^{p, q}}$,
\[
g \star L = g \circ L \circ g^{-1} = \sum_{s \in \CS} \widehat{N^{s,s}} + t \cdot \widehat{Z^{p, q}} = N + t \cdot \widehat{Z^{p, q}}
\]
as elements of $\rep_Q^\dd(\Bbbk[t,t^{-1}])$.
Hence the claim follows from the fact that $\ov{\CO}_M$ is a $\GL_\dd$-invariant subscheme of $\rep_Q^\dd$
and by Lemma~\ref{propertyclosedsub} applied to the canonical injective homomorphism $\Bbbk[t]\to\Bbbk[t,t^{-1}]$.

Applying the homomorphism $\Bbbk[t]\to\Bbbk[\varepsilon]$ sending $t$ to $\varepsilon$ we get that
$N+\varepsilon\cdot\widehat{Z^{p,q}}$ belongs to $\ov{\CO}_M(\Bbbk[\varepsilon])$.
\end{proof}

The above lemma gives a method of detecting vectors tangent to $\ov{\CO}_M$.
This method is sufficient for the representations of the Dynkin quivers of type $\BA$ as the proposition below shows.
Obviously this proposition follows also immediately from Theorem~\ref{mainforA}.

\begin{prop} \label{prop An}
Let $Q$ be a Dynkin quiver of type $\BA$ and $  M\in\rep(Q)$.
Then $\ov{\CO}_M(\Bbbk[\varepsilon])=\CC_M(\Bbbk[\varepsilon])$.
In other words, $T_N\ov{\CO}_M=T_N\CC_M$ for any $N$ in $\ov{\CO}_M(\Bbbk)=\CC_M(\Bbbk)$.
\end{prop}

\begin{proof}
Let $N \in \ov{\CO}_M(\Bbbk)$ and fix a decomposition $N=\bigoplus N^s$ of $N$ such that each $N^s$ is indecomposable.

Choose $Z \in T_N\CC_M=\BZ^1_{M,N}(N,N)$.
Since $Z=\sum_{p,q\in\CS}\widehat{Z^{p,q}}$ (see~\eqref{eq_sum}),
it is sufficient to show that $\widehat{Z^{p,q}}\in T_N\ov{\CO}_M$, for all $p, q \in \CS$.
Fix such $p$ and $q$.
We may assume $\widehat{Z^{p,q}}\not\in\BB^1_Q(N,N)$, as $\BB^1_Q(N,N)=T_N\CO_N\subseteq T_N\ov{\CO}_M$.
Since $\BZ^1_{M,N}(?,?)$ and $\BB^1_Q(?,?)$ are subfunctors of $\BZ^1_Q(?,?)$,
$Z^{p,q}$ belongs $\BZ^1_{M,N}(N^q,N^p)$ but not to $\BB^1_Q(N^q,N^p)$.
In particular, $\Ext^1_Q(N^q,N^p)$ is non-zero.
Since $Q$ is a Dynkin quiver, $N^q$ is not isomorphic to $N^p$, hence $q\neq p$.

It follows from the definition of $\BZ^1_{M,N}(?,?)$ that
$\supp(\delta_{\sigma(N^p,Z^{p,q},N^q)})\subseteq\supp(\delta_{M,N})$.
Combining Lemmas~\ref{homifpath}(4) and~\ref{fundamental}(3) we get that the values
of the function $\delta_{\sigma(N^p,Z^{p,q},N^q)}$ do not exceed $1$.
Hence we conclude the inequality $\delta_{\sigma(N^p,Z^{p,q},N^q)}\leq\delta_{M,N}$,
thus $\widehat{Z^{p,q}}$ belongs to $T_N\ov{\CO}_M$, by Lemma~\ref{constrtangent}.
\end{proof}

The above method does not extend to the Dynkin quivers of types $\BD$ and $\BE$.
A reason for this is that for these quivers there exist short exact sequences
$\sigma$ with indecomposable end terms such that the functions $\delta_\sigma$
attain values larger than $1$.

\begin{prop} \label{finalmainprop}
Let $Q$ be a Dynkin quiver of type $\BD$ and $N=\bigoplus_{s\in\CS}N^s$ be a decomposition
of $N \in \ov{\CO}_M(\Bbbk)$ such that each representation $N^s$ is indecomposable.
Let $Z^{p,q}\in\BZ_{M,N}^1(N^q,N^p)$, for $p\neq q$ in $\CS$,
be such that the inequality $\delta_{\sigma(N^p,Z^{p,q},N^q)}\leq\delta_{M,N}$ does not hold.

Then there is an index $r$ in $\CS\setminus\{p,q\}$ and a homomorphism $h^{q,r}$ in $\Hom_Q(N^r,N^q)$
such that for $Y^{p,r}=Z^{p,q} \circ h^{q,r}$ the following conditions hold:
\begin{enumerate}
\item[\textup{(1)}] The point $L=N+\widehat{Y^{p,r}}$ in $\rep_Q^\dd(\Bbbk)$ belongs to $\ov{\CO}_M(\Bbbk)$.
\item[\textup{(2)}] $N$ is a proper degeneration of $L$, i.e.\ $\ov{\CO}_N\varsubsetneq\ov{\CO}_L$.
 In particular, $\dim\CO_N<\dim\CO_L$.
\item[\textup{(3)}] The point $L+\varepsilon\cdot\widehat{Z^{p,q}}$ in $\rep_Q^\dd(\Bbbk[\varepsilon])$
 belongs to $\CC_M(\Bbbk[\varepsilon])$
 (equivalently, $\widehat{Z^{p,q}}\in T_L\CC_M)$.
\end{enumerate}
\end{prop}

\begin{proof}
We apply Corollary~\ref{maincortriang} for $\sigma=\sigma(N^p,Z^{p,q},N^q)$.
Hence there is an index $r$ in $\CS$ and a homomorphism $h^{q,r}$ in $\Hom_Q(N^r,N^q)$
such that for $Y^{p,r}=Z^{p,q} \circ h^{q,r}$ the following conditions hold:
\begin{enumerate}
\item[(i)] The sequence $\sigma(N^p,Y^{p,r},N^r)\colon\;
0\to N^p\xrightarrow{\bsmatrix{1\\ 0}}\bsmatrix{N^p&Y^{p,r}\\ 0&N^r}\xrightarrow{\bsmatrix{0&1}}N^r\to 0$
does not split. Equivalently, $Y^{p,r}$ does not belong to $\BB^1_Q(N^r,N^p)$.
\item[(ii)] $\delta_{\sigma(N^p,Y^{p,r},N^r)}\leq\delta_{M,N}$.
\item[(iii)] $\supp(\delta_{\sigma(N^p,Z^{p,q},N^q)}-\delta_{\sigma(N^p,Y^{p,r},N^r)})
\subseteq\supp(\delta_{M,N}-\delta_{\sigma(N^p,Y^{p,r},N^r)})$.
\end{enumerate}

We claim that $r\neq p,q$.
Note that (i) means that $Y^{p,q} = Z^{p,q} \circ h^{q,r}$ does not belong to $\BB^1_Q(N^r,N^p)$.
In particular $\Ext^1_Q(N^r,N^p)$ is non-zero, hence $r\neq p$.
Moreover $h^{q,r}$ is non-zero.
Suppose that $q=r$.
By Lemma~\ref{only1nonzero}(2), $\End_Q(N^q)=\Bbbk$,
hence $h^{q,r}$ would be an isomorphism, thus $\delta_{\sigma(N^p,Y^{p,r},N^r)}=\delta_{\sigma(N^p,Z^{p,q},N^q)}$.
Then (ii) contradicts the assumptions on $\delta_{\sigma(N^p,Z^{p,q},N^q)}$.
Consequently, $r\neq q$, and the claim is proved.

Now (1) follows from (ii) and Lemma~\ref{constrtangent},
and (2) is a consequence of (i) and Lemma~\ref{deltasigmapq}.
In order to prove (3) it suffices by~\eqref{descriptionCMtwo} to show that
\begin{enumerate}
\item[(iv)] $\supp(\delta_{\sigma(L,\widehat{Z^{p,q}},L)})\subseteq\supp(\delta_{M,L})$.
\end{enumerate}

As in the proof of Lemma~\ref{deltasigmapq}, we see that
$L$ is isomorphic to $\bigoplus_{s\neq p,r}N^s\oplus\bsmatrix{N^p&Y^{p,r}\\ 0&N^r}$.
Observe that $\sigma(L,\widehat{Z^{p,q}},L)$ is the direct sum of the following split sequence
\[
0\to\bigoplus_{s\neq p,r}N^s\xrightarrow{\bsmatrix{1\\ 0\\ 0}}
\bigoplus_{s\neq p,r}N^s\oplus\bigoplus_{s\neq p,q,r}N^s\oplus\bsmatrix{N^p&Y^{p,r}\\ 0&N^r}
\xrightarrow{\bsmatrix{0&1&0\\ 0&0&1}}\bigoplus_{s\neq p,q,r}N^s\oplus\bsmatrix{N^p&Y^{p,r}\\ 0&N^r}\to 0
\]
and the sequence
\[
0\to\bsmatrix{N^p&Y^{p,r}\\ 0&N^r}\xrightarrow{\bsmatrix{1&0\\ 0&1\\ 0&0}}
\bsmatrix{N^p&Y^{p,r}&Z^{p,q}\\ 0&N^r&0\\ 0&0&N^q}\xrightarrow{\bsmatrix{0&0&1}}N^q\to 0.
\]
Observe that the map
\[
\bsmatrix{1^p&0&0\\ 0&h^{q,r}&1^q\\ 0&1^r&0}\colon\bsmatrix{N^p&Z^{p,q} \circ h^{q,r}&Z^{p,q}\\ 0&N^r&0\\ 0&0&N^q}
\to\bsmatrix{N^p&Z^{p,q}\\0&N^q}\oplus N^r
\]
is an isomorphism of representations, hence
\[
\delta_{\sigma(L,\widehat{Z^{p,q}},L)}
=\delta_{N^r\oplus\bsmatrix{N^p&Z^{p,q}\\0&N^q},N^q\oplus\bsmatrix{N^p&Y^{p,r}\\ 0&N^r}}
=\delta_{\sigma(N^p,Z^{p,q},N^q)}-\delta_{\sigma(N^p,Y^{p,r},N^r)}.
\]
Combining this equality with (iii) and Lemma~\ref{deltasigmapq}(1) we get (iv), which finishes the proof.
\end{proof}

\begin{proof}[Proof of Theorem~\ref{main}]
Let $Q$ be a Dynkin quiver of type $\BD$.
By decreasing induction on $\dim \CO_N$ we prove that
$T_N\ov{\CO}_M=T_N\CC_M$ for all $N\in\ov{\CO}_M(\Bbbk)=\CC_M(\Bbbk)$.
Choose $N\in\ov{\CO}_M(\Bbbk)$ and
assume $T_L\ov{\CO}_M=T_L\CC_M$ for all $L\in\ov{\CO}_M$ with $\dim\CO_L>\dim\CO_N$.
We also fix a decomposition $N=\bigoplus N^s$ of $N$ such that each $N^s$ is indecomposable.

Let $Z \in \BZ^1_{M,N}(N,N)$.
Similarly as in the proof of Proposition~\ref{prop An} it is enough to show that $\widehat{Z^{p, q}} \in T_N\ov{\CO}_M$, for all $p, q \in \CS$.
By repeating arguments from this proof, we may assume $p \neq q$ and
the inequality $\delta_{\sigma(N^p,Z^{p,q},N^q)}\leq\delta_{M,N}$ does not hold,
thus we may apply Proposition~\ref{finalmainprop}.
In particular, there exists $Y^{p,r} \in \BZ_Q^1(N^r,N^p)$ such that
$L=N+\widehat{Y^{p,r}}$ belongs to $\ov{\CO}_M(\Bbbk)$, $\dim\CO_L>\dim\CO_N$, and $\widehat{Z^{p,q}}$ belongs to $T_L\CC_M$.
But the latter equals $T_L\ov{\CO}_M$ by induction hypothesis, hence
$N+\widehat{Y^{p,r}}+\varepsilon\cdot\widehat{Z^{p,q}}$ belongs to $\ov{\CO}_M(\Bbbk[\varepsilon])$.

We claim that the point $N+t\cdot\widehat{Y^{p,r}}+\varepsilon\cdot\widehat{Z^{p,q}}$
of $\rep_Q^\dd(\Bbbk[t]\otimes\Bbbk[\varepsilon])$
belongs to $\ov{\CO}_M(\Bbbk[t]\otimes\Bbbk[\varepsilon])$.
Consider the element $g$ in $\GL_\dd(\Bbbk[t,t^{-1}]\otimes\Bbbk[\varepsilon])$
given by $g=\sum_{s\neq r} \widehat{1^s}+t^{-1} \cdot \widehat{1^r}$.
Then $g^{-1}=\sum_{s\neq r} \widehat{1^s}+t \cdot \widehat{1^r}$.
Moreover,
\begin{multline*}
g \star (N+\widehat{Y^{p,r}}+\varepsilon\cdot\widehat{Z^{p,q}}) =
g \circ \left( \sum_{s\in \CS} \widehat{N^{s,s}} +\widehat{Y^{p,r}}+\varepsilon\cdot\widehat{Z^{p,q}}\right) \circ g^{-1}
\\
= \sum_{s\in \CS} \widehat{N^{s,s}} + t \cdot \widehat{Y^{p,r}}+\varepsilon\cdot\widehat{Z^{p,q}}
= N+t\cdot\widehat{Y^{p,r}}+\varepsilon\cdot\widehat{Z^{p,q}}
\end{multline*}
as elements of $\rep_Q^\dd(\Bbbk[t,t^{-1}]\otimes\Bbbk[\varepsilon])$.
Now the claim follows from the fact that $\ov{\CO}_M$ is a $\GL_\dd$-invariant subscheme of $\rep_Q^\dd$
and by Lemma~\ref{propertyclosedsub} applied to the canonical injective homomorphism
$\Bbbk[t]\otimes\Bbbk[\varepsilon]\to\Bbbk[t,t^{-1}]\otimes\Bbbk[\varepsilon]$.

Applying the homomorphism $\Bbbk[t]\otimes\Bbbk[\varepsilon]\to\Bbbk[\varepsilon]$ sending
$1\otimes\varepsilon$ to $\varepsilon$ and $t\otimes 1$ to $0$ we get that
$N+\varepsilon\cdot\widehat{Z^{p,q}}$ belongs to $\ov{\CO}_M(\Bbbk[\varepsilon])$.
Consequently, $\widehat{Z^{p,q}}\in T_N\ov{\CO}_M$,
which finishes the proof.
\end{proof}

\bibsection

\bigskip

\noindent
Grzegorz Bobi\'nski and Grzegorz Zwara\\
Faculty of Mathematics and Computer Science\\
Nicolaus Copernicus University\\
Chopina 12/18, 87-100 Toru\'n, Poland\\
E-mail: gregbob@mat.umk.pl, gzwara@mat.umk.pl

\end{document}